\documentclass[a4paper,DIV=12,10pt]{scrartcl}

\usepackage[utf8]{inputenc}
\usepackage{amsfonts,amssymb,amsmath,amsthm}
\usepackage{epsfig}

\usepackage{graphics}
\usepackage{subfig}
\usepackage{color}

\usepackage{tikz}
\usetikzlibrary{plotmarks}
\usetikzlibrary{shapes,arrows}
\usetikzlibrary{positioning}
\usetikzlibrary{trees,mindmap}
\usetikzlibrary{fadings,shapes.arrows,shadows}
\usetikzlibrary{decorations.pathreplacing}
\usetikzlibrary{calc}

\usepackage{booktabs}
\usepackage{cite}

\newcommand{\QH}{Q_n^{\mathrm{H}}}
\newcommand{\QG}{Q_n^{\mathrm{G}}}
\newcommand{\QGL}{Q_n^{\mathrm{GL}}}
\newcommand{\IH}{I_{\tau}^{\mathrm{H}}}
\newcommand{\IG}{I_{\tau}^{\mathrm{G}}}
\newcommand{\IGL}{I_{\tau}^{\mathrm{GL}}}
\newcommand{\It}{I_{\tau}^{k+1}}
\newcommand{\Rt}{R_{\tau}^k}
\newcommand{\BGL}{B_n^{\mathrm{GL}}}
\newcommand{\cGr}{\mathrm{cG(\mathit{r})}}
\newcommand{\cGPone}{\mathrm{cGP\text{-}C^1(\mathit{k})}}
\newcommand{\cGPtwo}{\mathrm{cGP\text{-}C^2(\mathit{k})}}

\newcommand{\llangle}{\langle\!\langle}
\newcommand{\rrangle}{\rangle\!\rangle}

\newcommand{\vertiii}[1]{%
{\left\vert\kern-0.25ex\left\vert\kern-0.25ex\left\vert #1 
 \right\vert\kern-0.25ex\right\vert\kern-0.25ex\right\vert}}

\newcommand{\ud}{\,\mathrm{d}}
\newcommand{\R}{\mathbb{R}}
\newcommand{\N}{\mathbb{N}}

\numberwithin{equation}{section}
\numberwithin{figure}{section}
\numberwithin{table}{section}

\renewcommand{\P}{\mathbb{P}}
\newcommand{\ur}{\hat{u}}
\newcommand{\Urth}{\widehat{U}_{\tau,h}}
\newcommand{\urth}{\hat{u}_{\tau,h}}
\newcommand{\fr}{\hat{f}}
\newcommand{\frtau}{\hat{f}_{\tau}}
\newcommand{\Frtau}{\widehat{F}_{\tau}}
\newcommand{\dt}{\partial_t}

\newcommand{\Pin}{\Pi^{k-2}_n}

\newtheorem{defi}{Definition}[section]
\newtheorem{thm}[defi]{Theorem}
\newtheorem{lem}[defi]{Lemma}
\newtheorem{rem}[defi]{Remark}
\newtheorem{cor}[defi]{Corollary}
\newtheorem{ass}[defi]{Assumption}
\newtheorem{prob}[defi]{Problem}

\makeatletter
\let\@fnsymbol\@arabic
\makeatother

\begin{document}

\title{Galerkin--collocation approximation in time for the wave equation and its post-processing}

\author{M.\ Anselmann$^\ast$, M.\ Bause$^{\ast,}$\thanks{bause@hsu-hh.de (corresponding 
author)}\;, S.\ Becher$^\dag$, G.\ Matthies$^\dag$\\
{\small $^\ast$ Helmut Schmidt University, Faculty of 
Mechanical Engineering, Holstenhofweg 85,}\\ 
{\small 22043 Hamburg, Germany}\\
{\small $^\dag$ Technical University of Dresden, Institute of Numerical
Mathematics,}\\ 
{\small 01062 Dresden, Germany}
}


\maketitle

\begin{abstract}
\textbf{Abstract.}
We introduce and analyze a class of Galerkin--collocation discretization schemes in time for the
wave equation. Its conceptual basis is the establishment of a direct connection between the
Galerkin method for the time discretization and the classical collocation methods, with the
perspective of achieving  the accuracy of the former with reduced computational costs provided
by the latter in  terms of less complex linear algebraic systems. Continuously differentiable in
time discrete solutions are obtained by the application of a special quadrature rule involving
derivatives. Optimal order error estimates are proved for fully discrete approximations based on
the Galerkin--collocation approach. Further, the concept of Galerkin--collocation approximation
is extended to twice continuously differentiable in time discrete solutions. A direct connection
between the two families by a computationally cheap post-processing is
presented. The error
estimates are illustrated by numerical experiments.
\end{abstract}



\section{Introduction}
\label{Sec:Introduction}

In this work we introduce and analyze a Galerkin--collocation (cGP--C$^k$,
$k\in\{1,2\}$) approach in time combined with a continuous Galerkin (cG)
finite element method in space to  approximate the solution to the second
order hyperbolic wave problem
\begin{equation}
\label{Eq:IBVP}
\begin{array}{r@{\;}c@{\;}l@{\hspace*{2ex}}l}
\dt^2 u- \Delta u & = & f & \mbox{in } \, \Omega \times (0,T]\,,\\[1ex]
u & = & 0 & \mbox{on } \; \partial \Omega \times (0,T]\,,\\[1ex]
u(\cdot ,0) = u_0\,, \;\; \dt u (\cdot ,0) & = & u_1 & \mbox{in } \, 
\Omega\,,
\end{array}
\end{equation}

with  C$^k$ regular functions in time. In \eqref{Eq:IBVP}, $T>0$ denotes
some final time and $\Omega$ is a polygonal or polyhedral bounded domain in
$\R^d$, with $d=2$ or $d=3$. The function $f: \Omega \times (0,T] \to \R$
and the initial values $u_0,u_1: \Omega \to \R$ are given data. The system
\eqref{Eq:IBVP} is studied as a prototype model for more sophisticated wave
phenomena of practical interest like, for instance, elastic wave
propagation governed by the Lam\'e--Navier equations, the Maxwell system,
or wave equations in coupled systems such as fluid-structure interaction
and fully dynamic poroelasticity \cite{MW12}. 

Our modification of the standard continuous Galerkin--Petrov method (cGP)
for time discretization (cf., e.g.,  \cite{AM89,BKRS18,BL94,KM04}) and the
innovation of this work comes through imposing collocation conditions
involving the discrete solution's derivatives at the discrete time nodes
while on the other hand downsizing the test space of the discrete
variational problem compared with the standard cGP approach. This idea was
recently introduced in \cite{BMW17} by two of the authors of this work for
first-order systems of ordinary differential equations. We refer to our
schemes as Galerkin--collocation methods. The collocation equations at the
discrete time nodes then enable us to ensure regularity of higher order in
time of the discrete solutions. A further key ingredient in the
construction of the Galerkin--collocation approach comes through the
application of a special quadrature formula, investigated in~\cite{JB09},
and the definition of a related interpolation operator for the right-hand
side term of the variational equation. Both of them use derivatives of the
given function. The Galerkin--collocation schemes rely in an essential way
on the perfectly matching set of polynomial spaces (trial and test space),
quadrature formula, and interpolation operator. For the discretization of
the spatial variables a continuous finite element approach is used here.
This is done for the sake of brevity. Usually, discontinuous Galerkin
methods are preferred; cf.\ \cite{AB19,B08,K14}. Beyond the higher order
regularity in the time, the Galerkin--collocation schemes offer appreciable
advantages for the solution of the arising linear systems by a favorable
impact on the matrix block structure; cf.\ \cite{AB19} for details.

For the subclass of discrete solutions being once continuously
differentiable in time an error analysis with optimal order error estimates
in time and space and in various norms is given. We will stress the key ideas of our 
error analysis and present a fundamental concept for analyzing generalized Galerkin
approximations to wave problems. One key point of our convergence proof for second-order 
hyperbolic problems is the weak stability result of Lemma~\ref{Lem:Stab_2}. Compared with 
usual stability results for parabolic problems or for first-order hyperbolic
problems (cf., e.g., \cite[Lemma 4.2]{ES16}) a stability is obtained such
that in the resulting error analysis some contributions can no longer be
absorbed by terms on the left-hand side of the error inequality like it is
typically done. Therefore, to prove error estimates of optimal order, the
error in the time derivatives $(\dt u^0_{\tau,h}, \dt u^1_{\tau,h})$ for
the discrete approximation pair $(u^0_{\tau,h},u^1_{\tau,h})$ of $(u,\dt
u)$ is bounded firstly. For this, a variational problem that is satisfied
by $(\dt u^0_{\tau,h},\dt u^1_{\tau,h})$ is identified. Then, a minor
extension of a result of \cite{KM04} becomes applicable to the thus
obtained problem. This yields an estimate for $\dt u - \dt u^0_{\tau,h}$
and $\dt^2 u - \dt u^1_{\tau,h}$. These auxiliary results then enable us to
prove the desired optimal-order error estimates for $u-u_{\tau,h}^0$ and
$\dt u- u_{\tau,h}^1$.

Space-time finite element methods with continuous and discontinuous
discretizations of the time and space variables for parabolic and
hyperbolic problems are well-known and have been studied carefully in the
literature; cf., e.g.,
\cite{ABM17,AMTX11,AM89,BL94,BG10,BKRS18,BRK17_2,BRK17,CF12,DG14,DFW16,F16,FP96,GN16,
HST11,HST13,J93,KM04,K14,LMN16,S15,T06} and the references therein. The
space-time approaches of these works differ by the choices of the trial
and, in particular, of the test spaces. Depending on the construction of
the test basis functions, either time-marching schemes defined by local
problems on the respective subintervals $(t_{n_-1},t_n]$ of $(0,T]$ (cf.,
e.g., \cite{ABM17,AMTX11,BKRS18,HST11,HST13,DFW16}) or schemes where all
time steps are solved simultaneously (cf., e.g., \cite{DFW16,GN16,SH18})
are obtained. Here, by choosing basis test functions supported on a single
subinterval $(t_{n_-1},t_n]$, we end up with a time-marching approach.
Further, strong relations between cGP schemes, collocation, and
Runge--Kutta methods have been observed. In \cite{AMN09,AMN11} they are
studied thoroughly. Moreover, nodal superconvergence properties of the cGP
method are known; cf.\ \cite[Eq.\ (2.2)]{AMN11}. In a recent work
\cite{BKRS18}, co-authored by one author of this work, a recursive
post-processing of the original continuous in time cGP solution is
presented and analyzed. The post-processed approximation is built on each
time interval upon the Gauss--Lobatto quadrature points of the actual time
interval, at which the classical cGP solution is superconvergent with one
extra order of accuracy. On the one hand, the post-processing lifts the
superconvergence of the original cGP solution at the Gauss--Lobatto
quadrature points to all points of the time interval by adding a higher
order correction term which vanishes at the Gauss--Lobatto quadrature
points. On the other hand, the post-processing, which is done sequentially
on the advancing time intervals and is of low computational costs, yields a
numerical approximation that is globally $C^1$-regular in time. In
\cite[Subsec.\ 3.2]{ES16} and \cite[p.\ 494]{MN06}, similar post-processing
techniques and lifting operators were studied for discontinuous Galerkin
approximations in time. The post-processing can nicely be exploited, for
instance, for an a-posteriori error control in time and an adaptive choice
of the time mesh. We explicitly note that in contrast to \cite{BKRS18},
where the continuous differentiability is obtained by a post-processing of
the continuous Galerkin--Petrov approximation, the higher order regularity
in time that is built in this work is an inherent part of the construction
of the discrete solution itself. This demands a different quadrature
formula and interpolation operator for the right-hand side function.

This work is organized as follows. In Sec.~\ref{Sec:NotPrem} we introduce
our notation and summarize preliminaries. In particular, quadrature
formulas and related interpolation operators are introduced. In
Sec.~\ref{Sec:CG} our  class of Galerkin--collocation schemes is presented.
In Sec.~\ref{Sec:PrepAna} some auxiliary results for our error analysis are
provided. Sec.~\ref{Sec:ErrAna} contains our error analysis for our family
of once continuously differentiable in time Galerkin--collocation methods.
In Sec.~\ref{Sec:EngCons} the conservation of energy by the numerical
schemes is studied. In Sec.\ \ref{Sec:PP} our construction principle is
extended to define a class of twice continuously differentiable in time
Galerkin-collocation approximation schemes for the wave equation. A link to
the first class of schemes by a post-processing procedure is presented.
Finally, in Sec.~\ref{Sec:NumExp} our error estimates are illustrated and
confirmed by numerical experiments.

\section{Notation and preliminaries}
\label{Sec:NotPrem}

\subsection{Function spaces and evolution form of continuous problem}

We use standard notation. $H^m(\Omega)$ is the Sobolev space of $L^2(\Omega)$ 
functions with derivatives up to order $m$ in 
$L^2(\Omega)$ and $\langle \cdot,\cdot \rangle$ denotes the inner product in $L^2(\Omega)$. 
Further, $\llangle \cdot, \cdot \rrangle$ defines the $L^2$ inner product 
on the product space $L^2(\Omega)\times L^2(\Omega)$. We let $H^1_0(\Omega):=\{u\in 
H^1(\Omega) \::\: u=0 \mbox{ on } \partial \Omega\}$. For short, we put 
\[
H:=L^2(\Omega)\qquad \text{and} \qquad V:= H^1_0(\Omega)\,.
\]
We denote by $V'$ the dual space of $V$ and use the notation
\begin{align*}
\| \cdot \| := \| \cdot\|_{L^2(\Omega)}\,,\qquad 
\| \cdot \|_m := \| \cdot\|_{H^m(\Omega)}, \quad m \in \N,
\end{align*} 
for the norms of the Sobolev spaces where we do not differ between the scalar- 
and vector-valued cases. Throughout, the meaning will be obvious from the context.
For a Banach space $B$, we let $L^2(0,T;B)$, $C([0,T];B)$, and $C^m([0,T];B)$, $m\in\N$, be
the Bochner spaces of $B$-valued functions, equipped with their natural norms. 
For a subinterval $J\subseteq [0,T]$, we use the notations
$L^2(J;B)$, $C^m(J;B)$, and $C^0(J;B):= C(J;B)$.

In what follows, for non-negative numbers $a$ and $b$, the expression 
$a\lesssim b$ stands for the inequality $a \leq C\, b$ with a generic constant $C$ that 
is independent of the sizes of the spatial and temporal meshes. The value of $C$ can depend on 
the regularity of the space mesh, the polynomial degrees used for the space-time 
discretization, and the data (including $\Omega$). 

For any given $u\in V$, let the operator $A:V\to V'$ be uniquely defined by
\[
\langle Au,v\rangle := \langle \nabla u, \nabla v\rangle \qquad \forall v\in V\,,
\]
where $\langle\cdot,\cdot\rangle$ on the left-hand side is understood as duality pairing
between $V'$ and $V$.
Further, we denote by $\mathcal{A}: V\times H\to H\times V'$ the operator
\[
\mathcal{A} = \begin{pmatrix}
0 & -I \\ A & 0 
\end{pmatrix}
\]
with the identity mapping $I: H\to H$. We let 
\begin{equation*}
X:=L^2(0,T;V)\times L^2(0,T;H)\,. 
\end{equation*}
Introducing the unknowns $u^0=u$ and $u^1=\dt u$, problem \eqref{Eq:IBVP} can be 
recovered in evolution form.

\begin{prob}
\label{Prob:0}
\mbox{}\\
Let $f \in L^2(0,T;H)$ and $(u_0,u_1)\in V\times H$ be given and $F=(0,f)$. Find 
$U=(u^0,u^1)\in X$ 
such that
\begin{equation}
\label{Eq:EP1}
 \dt U + \mathcal{A} U = F \quad \mathrm{in}\;\; (0,T)\,, \quad U(0) = U_0=
 (u_0,u_1)\,.
 \end{equation}
\end{prob}

Problem \eqref{Prob:0} admits a unique solution $U\in X$ and the 
mapping $(f,u_0,u_1) \mapsto \left(u^0,u^1\right)$ is a linear continuous map from
$L^2(0,T;H)\times V \times H$ to $X$; cf.\ \cite[p.\ 273, Thm.\ 1.1]{L71}.
Further, $u^0 \in C([0,T];V)$ and $u^1 \in C([0,T];H)$ are 
satisfied; cf.\ \cite[p.\ 275, Thm.\ 8.2]{LM72}. It follows from~\eqref{Eq:EP1} that
$\dt u^1 \in L^2(0,T;V')$.

\begin{ass}
\label{Ass:DataSol}
\mbox{}\\
Throughout, we tacitly assume that the solution $u$ of \eqref{Eq:IBVP} satisfies all 
the additional regularity conditions that are required in our analysis. In addition, let
$f\in C^s([0,T];H)$ for some sufficiently large parameter $s\in \N$ be satisfied. 
\end{ass}

The first of the conditions in Assumption \ref{Ass:DataSol} implies further assumptions 
on the data $f,u_0,u_1$ and the boundary $\partial \Omega$ of $\Omega$. Improved 
regularity results for solutions to the wave problem \eqref{Eq:IBVP} can be found in, 
e.g., \cite[Sec.\ 7.2]{E10}. The second condition in Assumption~\ref{Ass:DataSol} 
will allow us to apply an interpolation in time that is based on derivatives of the 
right-hand side function $f$.

\subsection{Time and space discretization}

For the time discretization, we decompose the time interval $I=(0,T]$ into $N$ subintervals 
$I_n=(t_{n-1},t_n]$, $n=1,\ldots,N$, where $0=t_0<t_1< \cdots < t_{N-1} < t_N 
= T$ such that $I=\bigcup_{n=1}^N I_n$. We put $\tau = \max_{n=1,\ldots, N} \tau_n$ with 
$\tau_n = t_n-t_{n-1}$. Further, the set $\mathcal{M}_\tau := \{I_1,\ldots, I_N\}$
of time intervals is called the time mesh. For a Banach space $B$ and any $k\in \N_0$, 
we let 
\begin{equation}
\label{Def:Pk}
\P_k(I_n;B) = \bigg\{w_\tau : I_n \to B \::\: w_\tau(t) = \sum_{j=0}^k 
W^j t^j \; \forall t\in I_n\,, \; W^j \in B\; \forall j \bigg\}\,.
\end{equation}
For an integer $k\in \N$, we introduce the space
\begin{equation}
\label{Eq:DefXk} 
X_\tau^k (B) := \left\{w_\tau \in C(\overline{I};B) \::\: w_\tau|_{I_n} \in
\P_k(I_n;B)\; \forall I_n\in \mathcal{M}_\tau \right\}
\end{equation}
of globally continuous functions in time 
and for an integer $l\in \N_0$ the space
\begin{equation*}
\label{Eq:DefYk}
Y_\tau^{l} (B) := \left\{w_\tau \in L^2(I;B) \::\: w_\tau|_{I_n} \in
\P_{l}(I_n;B)\; \forall I_n\in \mathcal{M}_\tau \right\}
\end{equation*}
of global $L^2$-functions in time.

For any non-negative integer $s$ and a function $w:I\to B$ that is piecewise
sufficiently smooth with respect to the time mesh $\mathcal{M}_{\tau}$, we define by
\begin{equation}
\label{Eq:Defw_In_bdr}
\dt^s w(t_n^+) := \lim_{t\to t_n+0} \dt^s w(t)
\qquad\text{and}\qquad
\dt^s w(t_n^-) := \lim_{t\to t_n-0} \dt^s w(t)
\end{equation}
the one-sided limits of the $s$th derivative of $w$.

For the space discretization, let $\mathcal{T}_h$ be a shape-regular mesh of $\Omega$ 
consisting of quadrilateral or hexahedral elements with mesh size $h>0$. For some 
integer $r\in \N$, let $V_h=V_h^{(r)}$ be the scalar finite element space given by 
\begin{equation}
\label{Eq:DefVh}
 V_h = V_h^{(r)}=\left\{v_h \in C(\overline{\Omega}) \::\: v_h{}|_T\in \mathbb{Q}_r(K) \, \forall K 
\in \mathcal{T}_h \right\}\cap H^1_0(\Omega)
\end{equation}

where $\mathbb{Q}_r(K)$ is the space defined by the multilinear reference mapping of
polynomials on  the reference element with maximum degree $r$ in each variable. Our
restriction in this work to continuous finite elements in space is only done for
simplicity and in order to  reduce the technical methodology of analyzing our
Galerkin--collocation discretization scheme to its key points. In the literature it has
been mentioned that discontinuous finite element methods in space offer appreciable
advantages over continuous ones for the discretization of wave equations; cf., e.g.,
\cite{A06,B08,K14,K15} and the references therein.

We denote by $P_h: H\to V_h$ the $L^2$-orthogonal projection onto $V_h$ such that for $w\in H$,
\begin{equation*}
 \langle P_h w, v_h \rangle = \langle w, v_h\rangle 
\end{equation*}
for all $v_h\in V_h$. The operator $R_h: V\to V_h$ defines the elliptic 
projection onto $V_h$ such that for $w\in V$, 
\begin{equation}
\label{Def:EllipProj}
 \langle \nabla R_h w, \nabla v_h \rangle = \langle \nabla w, \nabla v_h\rangle 
\end{equation}
for all $v_h\in V_h$. Finally, by $\mathcal{P}_h : H \times H\to V_h\times V_h$ we 
denote the $L^2$-projection onto the product space $V_h\times V_h$ and by
$\mathcal{R}_h : V \times V\to V_h\times V_h$ the elliptic projection onto the product 
space $V_h\times V_h$. Let $A_h: H^1_0(\Omega) \to V_h$ be the operator that is defined by 
\begin{equation}
\label{Eq:DefAh}
 \langle A_h w , v_h \rangle = \langle \nabla w, \nabla v_h\rangle 
\end{equation}
for all $v_h\in V_h$. Then, for $w \in V\cap H^2(\Omega)$ it holds that 
\begin{equation*}
 \langle A_h w , v_h \rangle = \langle \nabla w, \nabla v_h \rangle = \langle Aw 
,v_h\rangle 
\end{equation*}
for all $v_h\in V_h$. Hence, we have $A_h w = P_h Aw $ for $w\in V\cap H^2(\Omega)$. Let 
$\mathcal{A}_h: V\times H \to V_h \times V_h$ be 
defined by 
\begin{equation*}
\mathcal{A}_h = \begin{pmatrix}
0 & -P_h\\ A_h & 0 
\end{pmatrix}\,.
\end{equation*}
Hence, we have for $W=(w^0,w^1)\in \big(V\cap H^2(\Omega)\big)\times H$ that
\begin{equation*}
\llangle \mathcal{A}_h W, \Phi_h \rrangle = \langle -w^1 , \phi_h^0\rangle + \langle 
\nabla w^0 , \nabla \phi_h^1\rangle = \langle - w^1 , \phi_h^0\rangle + \langle A 
w^0, \phi_h^1\rangle = \llangle \mathcal{A} W, \Phi_h \rrangle
\end{equation*}
for all $\Phi_h = (\phi_h^0,\phi_h^1) \in V_h\times V_h$. This provides the consistency
\begin{equation}
\label{Eq:ConsistA_h}
\mathcal{A}_h W = \mathcal{P}_h \mathcal{A}W
\end{equation} 
of $\mathcal{A}_h$ on $\big(V\cap H^2(\Omega)\big)\times H$.

Finally, let $U_{0,h} \in V_h^2$ denote a suitable approximation of the
initial value $U_0 \in V \times H$ in~\eqref{Eq:EP1} that will we used as
the initial value $U_{\tau,h}(0)$ of the discrete solution. Further
restrictions will be made below. 

\subsection{Quadrature formulas and interpolation operators}
Throughout this work, the polynomial degree $k\ge 3$ is assumed to be fixed. Let
$\hat{t}^{\,\mathrm{H}}_1=-1$, $\hat{t}^{\,\mathrm{H}}_{k-1}=1$, and $\hat{t}^{\,\mathrm{H}}_s$,
$s=2,\dots,k-2$, be the roots of the Jacobi polynomial on $\widehat{I}:=[-1,1]$ with
degree $k-3$ associated to the weighting function $(1-\hat{t})^2(1+\hat{t})^2$.
Let $\widehat{I}^{\,\mathrm{H}}:C^1\big(\widehat{I};B\big)\to \P_k\big(\widehat{I};B\big)$
denote the Hermite interpolation operator with respect to point value and first
derivative at both $-1$ and $1$ as well as the point values at $\hat{t}_s^{\,\mathrm{H}}$,
$s=2,\dots,k-2$. By
\begin{equation}
\widehat{Q}^{\mathrm{H}}(\hat{g}) := \int_{-1}^{1}
\widehat{I}^{\,\mathrm{H}}(\hat{g})(\hat{t})\ud\hat{t}
\end{equation}
we define an Hermite-type quadrature on $[-1,1]$ which can be written as
\begin{equation}
\widehat{Q}^\mathrm{H}(\hat{g}) = \widehat{\omega}_L \hat{g}'(-1) + \sum_{s=1}^{k-1}
\widehat{\omega}_s \hat{g}(\hat{t}^{\,\mathrm{H}}_s) + \widehat{\omega}_R \hat{g}'(1)\,,
\end{equation}
where all weights are non-zero. Using the affine mapping $T_n:\widehat{I}\to \overline{I}_n$
with $T_n(-1) = t_{n-1}$ and $T_n(1) = t_n$, we obtain
\begin{equation}
\label{Eq:GLHF}
\QH(g) = \left(\frac{\tau_n}{2}\right)^2 \widehat{\omega}_L \dt g(t_{n-1}^+)
+ \frac{\tau_n}{2} \sum_{s=1}^{k-1} \widehat{\omega}_s g(t^\mathrm{H}_{n,s})
+ \left(\frac{\tau_n}{2}\right)^2 \widehat{\omega}_R \dt g(t_n^-)
\end{equation}
as Hermite-type quadrature formula on $I_n$, where $t^\mathrm{H}_{n,s} :=
T_n(\hat{t}^\mathrm{H}_s)$, $s=1,\dots,k-1$. We note that $\QH$ given in~\eqref{Eq:GLHF} integrates all polynomials up to degree $2k-3$ exactly, cf.~\cite{JB09}. Using
$\widehat{I}^{\,\mathrm{H}}$ and $T_n$, the local Hermite interpolation on $I_n$ is given by
\[
I_n^\mathrm{H} : C^1(\overline{I}_n;B)\to \P_k(\overline{I}_n;B) \,,
\qquad 
v\mapsto \big(\widehat{I}^{\,\mathrm{H}}(v\circ T_n)\big)\circ T_n^{-1}\,.
\]
Moreover, we define the global Hermite
interpolation $\IH:C^1(\overline{I};B)\to X_\tau^k(B)$ by means of 
\begin{equation}
\label{Eq:DefHIntOp}
\IH w|_{I_n} := I_n^{\mathrm{H}}(w|_{I_n})
\end{equation}
for all $n=1,\dots,N$.

In addition to Hermite-type interpolation and quadrature formula, Gauss and
Gauss--Lobatto quadrature formulas will be used. To this end, we denote by $\hat{t}_s^{\, \mathrm{G}}$,
$s=1,\dots,k-1$, the roots of the Legendre polynomial with degree $k-1$ and by
$\hat{t}_s^{\,\mathrm{GL}}$, $s=2,\dots,k-1$, the roots of the Jacobi polynomial on
$\widehat{I}$ with degree $k-2$ associated to the weighting function
$(1-\hat{t})(1+\hat{t})$. Furthermore, we set $\hat{t}_1^{\,\mathrm{GL}}=-1$ and
$\hat{t}_k^{\,\mathrm{GL}}=1$. The operators 
$\widehat{I}^{\,\mathrm{G}}:C(\widehat{I};B)\to P_{k-2}(\widehat{I};B)$ and 
$\widehat{I}^{\,\mathrm{GL}}:C(\widehat{I};B)\to P_{k-1}(\widehat{I};B)$ are the Lagrange
interpolation using the Gauss points $\hat{t}_s^{\,\mathrm{G}}$, $s=1,\dots,k-1$, and the
Gauss--Lobatto points $\hat{t}_s^{\,\mathrm{GL}}$, $s=1,\dots,k$, respectively.
We define by
\begin{equation}
\widehat{Q}^{\,\mathrm{G}}(\hat{g}) := \int_{-1}^1
\widehat{I}^{\,\mathrm{G}}(\hat{g})(\hat{t})\ud\hat{t}
\qquad\text{and}\qquad
\widehat{Q}^{\,\mathrm{GL}}(\hat{g}) := \int_{-1}^1
\widehat{I}^{\,\mathrm{GL}}(\hat{g})(\hat{t})\ud\hat{t}
\end{equation}
Gauss and Gauss--Lobatto quadrature formulas on $[-1,1]$ which are transformed to
\begin{equation}
\label{Eq:GF+GLF}
\QG(g) = \frac{\tau_n}{2} \sum_{s=1}^{k-1} \widehat{\omega}_s^\mathrm{G}
g(t_{n,s}^\mathrm{G})
\qquad\text{and}\qquad
\QGL(g) = \frac{\tau_n}{2} \sum_{s=1}^k \widehat{\omega}_s^\mathrm{GL}
g(t_{n,s}^\mathrm{GL})
\end{equation}
on $I_n$ by using the affine mapping $T_n$. The Gauss and Gauss--Lobatto formulas also integrate
polynomials up to degree $2k-3$ exactly. Local Lagrange-type interpolation operators on $I_n$
are given by
\begin{alignat*}{2}
I_n^\mathrm{G} : C(\overline{I}_n;B) &\to \P_{k-2}(\overline{I}_n;B)\,,
& \qquad 
v &\mapsto \big(\widehat{I}^{\,\mathrm{G}}(v\circ T_n)\big)\circ T_n^{-1}\,,\\
I_n^\mathrm{GL} : C(\overline{I}_n;B) &\to \P_{k-1}(\overline{I}_n;B)\,,
& \qquad 
v &\mapsto \big(\widehat{I}^{\,\mathrm{GL}}(v\circ T_n)\big)\circ T_n^{-1}\,.
\end{alignat*}
Furthermore, we define the global Lagrange interpolation operators
$\IG:C(\overline{I};B)\to Y_\tau^{k-2}(B)$ and $\IGL:C(\overline{I};B)\to X_\tau^{k-1}(B)$ by
\begin{equation*}
\IG w|_{I_n} := I_n^{\mathrm{G}}(w|_{I_n})
\qquad\text{and}\qquad
\IGL w|_{I_n} := I_n^{\mathrm{GL}}(w|_{I_n})
\end{equation*}
for all $n=1,\dots,N$.

\section{Galerkin--collocation discretization and 
auxiliaries}
\label{Sec:CG}

In this section we introduce the approximation of the wave problem~\eqref{Eq:EP1} by our
Galerkin--collocation approach that combines collocation conditions at the endpoint
$t_{n-1}$ and $t_n$ of the subintervals $I_n$ with variational equations for reduced test
spaces compared with the standard continuous finite element approximation of the wave
equation (cf.\ \cite{BKRS18,FP96,KM04}). A family of discrete solutions that are once 
continuously differentiable in time is obtained. For this family an
optimal order error analysis is then developed in Sec.~\ref{Sec:ErrAna}. For the sake of 
completeness and in order to show the impact of the collocation conditions, the standard 
continuous Galerkin approximation (cf.\ \cite{FP96,KM04}) of the wave 
problem~\eqref{Eq:EP1} is briefly recalled in Subsec.~\ref{Sec:cGP}.

\subsection{Space-time discretization with continuous Galerkin--Petrov method 
$\boldsymbol{\mathrm{cGP(\mathit{k})}}$}
\label{Sec:cGP}

For completeness and comparison, we briefly present the standard continuous
Galerkin--Petrov method of order $k\geq 1$ (in short, $\mathrm{cGP(\mathit{k})}$) as time discretization
applied to the evolution problem~\eqref{Eq:EP1}. For the space discretization, the
continuous Galerkin approach $\cGr$ in $V_h$, defined in 
\eqref{Eq:DefVh}, is used for the sake of simplicity. This yields the following fully 
discrete problem; cf., e.g., \cite{BKRS18,KM04} for details.

\begin{prob}[Global, fully discrete problem of $\mathrm{cGP(\mathit{k})}$--$\cGr$]
\label{Prob:SemiDis}
\mbox{}\\
Find $U_{\tau,h} \in \big(X_\tau^k (V_h)\big)^2$ such that $U_{\tau,h} (0) = 
U_{0,h}$ 
and
\begin{equation*}
\int_{0}^T \Big(\llangle \dt U_{\tau,h} , V_{\tau,h} \rrangle + 
\llangle \mathcal{A}_h U_{\tau,h} , V_{\tau,h} \rrangle \Big) \ud t = \int_{0}^T 
\llangle F,V_{\tau,h} \rrangle \ud t
\end{equation*}
for all $V_{\tau,h} \in \big(Y_\tau^{k-1} (V_h) \big)^2$.
\end{prob}

Both components of $U_{\tau,h}=(u_{\tau,h}^0,u_{\tau,h}^1)$ are computed in the 
same discrete space $X_\tau^k (V_h)$. By choosing test functions supported on a 
single subinterval $I_n$ and using the $(k+1)$-point Gauss--Lobatto quadrature formula,
we recast Problem~\ref{Prob:SemiDis} as a sequence of local 
problems on $I_n$.

\begin{prob}[Local, numerically integrated, fully discrete problem of
$\mathrm{cGP(\mathit{k})}$--$\cGr$ on $I_n$]
\label{Prob:SemiDisLoc}
\mbox{}\\
Find $U_{\tau,h}|_{I_n} \in \big((\P_k (I_n;V_h)\big)^2$
with $U_{\tau,h}(t_{n-1}^+)=U_{\tau,h}(t_{n-1}^-)$ for 
$n>1$ and $U_{\tau,h}(t_0^+)=U_{0,h}$ 
such that
\begin{equation*}
\label{Eq:FullDisLocal}
Q_{n,k+1}^{\mathrm{GL}} \big(\llangle \dt U_{\tau,h} , V_{\tau,h} \rrangle 
+ \llangle \mathcal{A}_h U_{\tau,h} , V_{\tau,h} \rrangle \big) = 
Q_{n,k+1}^{\mathrm{GL}} \big(\llangle F,V_{\tau,h}\rrangle\big)
\end{equation*}
for all $V_{\tau,h} \in (\P_{k-1} (I_n;V_h))^2$.
\end{prob}
In Problem~\ref{Eq:FullDisLocal} we use a Gauss--Lobatto quadrature formula with $k+1$ 
points, which is in contrast to $\QGL$ in \eqref{Eq:GF+GLF} that uses $k$ points. 
Furthermore, the quadrature formula on the left-hand side can be replaced by exact 
integration or by any quadrature formula which is exact for
polynomials of degree up to order $2k-1$.

\subsection{Space-time discretization with Galerkin--collocation method
$\boldsymbol{\cGPone}$}
\label{Subsec:cGPC1}

From now on we suppose that $k\ge 3$ is a fixed integer without always mentioning this 
explicitly.

\begin{prob}[Local, numerically integrated, fully discrete problem of $\cGPone$--$\cGr$
on $I_n$]
\label{Prob:DisLocalcGPC}
\mbox{}\\
Given $U_{\tau,h}(t_{n-1}^-)$ for $n>1$ and $U_{\tau,h}(t_0^-)=U_{0,h}$ for $n=1$,
find $U_{\tau,h}|_{I_n} \in \big(\P_k(I_n;V_h)\big)^2$ such that
\begin{subequations}
\label{eq:P33}
\begin{align}
\label{Eq:DisLocalcGPC_1}
U_{\tau,h}(t_{n-1}^+) & = U_{\tau,h}(t_{n-1}^-) \,,\\[1ex]
\label{Eq:DisLocalcGPC_2}
\dt U_{\tau,h}(t_{n-1}^+) & = - \mathcal{A}_h U_{\tau,h}(t_{n-1}^+)
+ \mathcal{P}_h F(t_{n-1}^+) \,, \\[1ex]
\label{Eq:DisLocalcGPC_3}
\dt U_{\tau,h}(t_{n}^-) & = - \mathcal{A}_h U_{\tau,h}(t_{n}^-) + \mathcal{P}_h F(t_{n}^-)\,,
\intertext{and}
\label{Eq:DisLocalcGPC_4}
\QH \Big(\llangle \dt U_{\tau,h} , V_{\tau,h} \rrangle 
& + \llangle \mathcal{A}_h U_{\tau,h} , V_{\tau,h} \rrangle \Big) = 
\QH \big(\llangle F,V_{\tau,h}\rrangle\big)
\end{align}
\end{subequations}
for all $V_{\tau,h} \in \big(\P_{k-3} (I_n;V_h)\big)^2$.
\end{prob}

For this scheme we make the following observations.

\begin{rem}
It directly follows from the definition of the scheme that $U_{\tau,h} \in 
\big(C^1(\overline{I}; V_h)\big)^2$ is satisfied. Instead of the condition 
\eqref{Eq:DisLocalcGPC_2} at $t_{n-1}^+$ we could also demand that
\begin{equation}
\label{Eq:PracCond:t_n-1} 
\dt U_{\tau,h}(t_{n-1}^+) = \dt U_{\tau,h}(t_{n-1}^-)\,,
\end{equation}
where we set $\dt U_{\tau,h}(t_0^-) = - \mathcal{A}_h U_{0,h} + \mathcal{P}_h F(0)$.

Since the time discretization is of Galerkin--Petrov type, we refer to it as a 
continuously differentiable Galerkin--Petrov approximation, for short $\cGPone$.

Compared to Problem~\ref{Prob:SemiDisLoc}, the test space of the variational
constraint~\eqref{Eq:DisLocalcGPC_4} reduces from $\big(\P_{k-1} (I_n;V_h)\big)^2$
to $\big(\P_{k-3} (I_n;V_h)\big)^2$. For $k=3$ the test space just becomes 
the set $\big(\P_{0} (I_n;V_h)\big)^2$ of piecewise constant functions in time. 
Introducing the collocation conditions~\eqref{Eq:DisLocalcGPC_2} and
\eqref{Eq:DisLocalcGPC_3} along with downsizing the test space of the 
variational condition impacts the block structure of the resulting linear 
algebraic system. By \eqref{Eq:PracCond:t_n-1} a condensation of internal degrees 
of freedom becomes feasible which leads to smaller algebraic systems and might simplify 
the future construction of efficient iterative solvers and preconditioners; 
cf.~\cite{AB19}. 
\end{rem}

The existence of a unique solution to Problem~\ref{Prob:DisLocalcGPC} can be proved
along the lines of \cite[p.\ 812, Thm.\ A.3]{BRK17} 
by using the equivalence of existence and uniqueness in the finite dimensional case. 

We state for the scheme~\eqref{eq:P33} the following auxiliary results. 

\begin{lem}
\label{Lem:DSVP}
The solution $U_{\tau,h}\in \big(X^k_{\tau,h}(V_h)\big)^2$ of
Problem~\ref{Prob:DisLocalcGPC} satisfies for $n=1,\ldots,N$ that
\begin{equation}
\label{Lem:DSVP_0}
\QGL\big(\llangle \dt U_{\tau,h},V_{\tau,h}\rrangle+ \llangle 
\mathcal{A}_h U_{\tau,h} , V_{\tau,h} \rrangle\big) = 
\QGL\big(\llangle I_\tau^{\mathrm H}F ,V_{\tau,h}\rrangle\big) 
\end{equation}
for all $V_{\tau,h}\in \big(\P_{k-2}(I_n;V_h)\big)^2$.
\end{lem}

We note that compared to Problem~\ref{Prob:DisLocalcGPC} the quadrature formula has 
been changed in Lemma \ref{Lem:DSVP}. In addition, the test space has been increased from 
$\P_{k-3}$ to $\P_{k-2}$.

\begin{proof}
For arbitrarily chosen $V_{\tau,h}\in \big(\P_{k-2}(I_n;V_h)\big)^2$,
there exists some $d_n = d_n(V_{\tau,h})\in V_h^2$ such that $V_{\tau,h}$ admits the
representation
\begin{equation*}
V_{\tau,h} = \widetilde V_{\tau,h} + d_n(V_{\tau,h})\psi_n
\end{equation*}
with 
\begin{equation*}
\widetilde V_{\tau,h} \in \big(\P_{k-3}(I_n;V_h)\big)^2
\quad \text{and}\quad 
\psi_n(t) = \prod_{\mu=2}^{k-1}(t-t_{n,\mu}^{\mathrm{GL}})\in \P_{k-2}(I_n)
\end{equation*}
where $t_{n,\mu}^{\mathrm{GL}}$, $\mu=2,\ldots,k-1$, denote the inner 
Gauss--Lobatto quadrature points on $\overline{I}_n$. From~\eqref{Eq:DisLocalcGPC_4} along 
with the exactness of the Hermite-type quadrature formula~\eqref{Eq:GLHF} for all 
polynomials in $\P_{2k-3}(I_n)$ and of the Gauss--Lobatto quadrature formula 
\eqref{Eq:GF+GLF} for all polynomials in $\P_{2k-3}(I_n)$, it follows that
\begin{equation}
\label{Eq:DSVP_1}
\begin{aligned}
\QGL \big(\llangle \dt U_{\tau,h} , \widetilde V_{\tau,h} 
\rrangle + \llangle \mathcal{A}_h U_{\tau,h} , \widetilde V_{\tau,h} \rrangle \big)
& = \QH\big(\llangle \dt U_{\tau,h} , \widetilde V_{\tau,h} 
\rrangle + \llangle \mathcal{A}_h U_{\tau,h} , \widetilde V_{\tau,h} \rrangle \big)\\[1ex]
& = \QH\big(\llangle F,\widetilde V_{\tau,h}\rrangle\big) = 
\QH\big(\llangle I_\tau^{\mathrm H}F,\widetilde V_{\tau,h}\rrangle\big)\\
& = \QGL\big(\llangle I_\tau^{\mathrm H}F,\widetilde V_{\tau,h}\rrangle\big)\,.
\end{aligned}
\end{equation}
Therefore, it remains to prove that 
\begin{equation}
\label{Eq:DSVP_2}
\QGL\big(\llangle \dt U_{\tau,h} , d_n \psi_n
\rrangle + \llangle \mathcal{A}_h U_{\tau,h} , d_n \psi_n \rrangle \big) = 
\QGL\big(\llangle I_\tau^{\mathrm H}F,d_n \psi_n\rrangle\big)
\end{equation}
is satisfied. Since $\psi_n$ vanishes in the interior Gauss--Lobatto quadrature nodes 
$t_{n,\mu}^{\mathrm{GL}}$, $\mu = 2,\ldots, k-1$, and the quantities
$\llangle \dt U_{\tau,h} , d_n \psi_n \rrangle + \llangle \mathcal{A}_h U_{\tau,h},
d_n \psi_n \rrangle$ and $\llangle I_\tau^{\mathrm H}F,d_n \psi_n\rrangle$ coincide
in the endpoints $t_{n-1}^{+}$ and $t_n^{-}$ by means of the
conditions~\eqref{Eq:DisLocalcGPC_2} and~\eqref{Eq:DisLocalcGPC_3}, the variational
problem~\eqref{Eq:DSVP_2} is satisfied.
Along with~\eqref{Eq:DSVP_1}, this proves the assertion~\eqref{Lem:DSVP_0}.
\end{proof}

Furthermore, the solution of Problem~\ref{Prob:DisLocalcGPC} fulfills an evolution problem on
$I$.

\begin{lem} 
\label{Lem:EvolProb}
The solution $U_{\tau,h}$ of Problem~\ref{Prob:DisLocalcGPC} satisfies
\begin{equation}
\label{Eq:EvolProb_0}
\dt U_{\tau,h} + \IGL\mathcal{A}_h U_{\tau,h} = \mathcal{P}_h 
\IGL \IH F 
\end{equation}
on the whole time interval $\overline{I}$.
\end{lem}

\begin{proof}
Since all quantities in~\eqref{Eq:EvolProb_0} are continuous on $\overline{I}$, it suffices
to prove the relation locally on each $I_n$. To this end, let $n\in \{1,\ldots,N\}$ be fixed.

From~\eqref{Eq:DisLocalcGPC_3} along with $t_{n,k}^\mathrm{GL}=t_n$ and the
interpolation properties of $\IGL$ and $\IH$, we get that
\begin{equation}
\label{Eq:EvolProb_3}
\dt U_{\tau,h}(t_{n,k}^\mathrm{GL}) + \IGL \mathcal{A}_h
U_{\tau,h}(t_{n,k}^\mathrm{GL})
- \mathcal{P}_h \IGL \IH F(t_{n,k}^\mathrm{GL}) = 0 \,.
\end{equation}
Using~\eqref{Lem:DSVP_0}, it follows that 
\begin{equation}
\label{Eq:EvolProb_1}
\QGL\big(\llangle \dt U_{\tau,h} + \IGL \mathcal{A}_h U_{\tau,h}
- \mathcal{P}_h \IGL \IH F, V_{\tau,h}\rrangle\big) \\[1ex]
= \QGL\big(\llangle \dt U_{\tau,h} + \mathcal{A}_h U_{\tau,h}
- \IH F, V_{\tau,h}\rrangle\big) = 0
\end{equation}
for all $V_{\tau,h}\in \big(\mathbb{P}_{k-2}(I_n;V_h)\big)^2$. Choosing now test
functions $V_{\tau,h}^i\in \big(\P_{k-2}(I_n;V_h)\big)^2$, $i=1,\ldots, k-1$,
such that $V_{\tau,h}^i(t_{n,\mu}^\mathrm{GL}) = \delta_{i,\mu}\Phi_h$, 
$\mu=1,\ldots,k-1$, and $\Phi_h\in V_h\times V_h$, the properties~\eqref{Eq:EvolProb_1}
and~\eqref{Eq:EvolProb_3} result in
\begin{equation}
\label{Eq:EvolProb_2}
\dt U_{\tau,h}(t_{n,i}^\mathrm{GL}) + \IGL \mathcal{A}_h
U_{\tau,h}(t_{n,i}^\mathrm{GL})
- \mathcal{P}_h \IGL \IH F(t_{n,i}^\mathrm{GL}) = 0 \,,
\quad \text{for } i = 1,\ldots ,k-1\,.
\end{equation}
Thus, by means of~\eqref{Eq:EvolProb_2} and~\eqref{Eq:EvolProb_3}, the polynomial
$\dt U_{\tau,h} + \IGL\mathcal{A}_h U_{\tau,h} - \mathcal{P}_h \IGL \IH F\in
\big(\P_{k-1}(I_n;V_h)\big)^2$ vanishes in the $k$ different nodes
$t_{n,i}^\mathrm{GL}$, $i = 1,\ldots,k$. Therefore, it 
vanishes for all $t\in \overline{I}_n$, which proves the local version
of~\eqref{Eq:EvolProb_0}. The statement of this lemma follows from the global
continuity.
\end{proof}

\begin{rem}
\label{Rem:rel12a}
The statements of Lemma~\ref{Lem:DSVP} and Lemma~\ref{Lem:EvolProb} are quite similar to
the statements of Lemma\ 4.4 and Lemma\ 3.11 given in~\cite{BKRS18}, respectively.
However, in contrast to the analysis of~\cite{BKRS18}, the pointwise
identity~\eqref{Eq:EvolProb_0} is not needed for the proof of~\eqref{Lem:DSVP_0} since the
collocation conditions~\eqref{Eq:DisLocalcGPC_2} and~\eqref{Eq:DisLocalcGPC_3} already
provide the needed additional information.
\end{rem}

\section{Preparation for the error analysis} 
\label{Sec:PrepAna}

We will use in our error analysis some interpolants in time introduced in
\cite{BKRS18,ES16}. To keep this work self-contained, their definition and some
auxiliaries are briefly summarized here. Remember that $k\geq 3$.

In the following, let $B$ be a Banach space satisfying $B\subset H$ and
$\ell\in \mathbb{N}$. We define for $n=1,\dots,N$ the local $L^2$-projections
$\Pi^{\ell}_n:L^2(I_n;B)\to
\P_{\ell}(I_n;B)$ by
\begin{equation}
\label{Def:Pi}
\int_{I_n} \langle \Pi^{\ell}_n w , q \rangle \ud t =
\int_{I_n} \langle w , q \rangle \ud t 
\qquad\forall\, q\in \P_{\ell}(I_n;B).
\end{equation}

Next, a special interpolant in time is constructed. To this end, we define
the Hermite interpolation operator $\It:C^1(\overline{I};B)\to
C^1(\overline{I};B)\cap X_{\tau}^{k+1}(B)$ by
\[
\It u(t_n) = u(t_n),\quad
\dt \It u(t_n) = \dt u(t_n)\,,
\qquad n=0,\dots,N,
\]
and
\[
\It u(t_{n,\mu}^\mathrm{GL}) = u(t_{n,\mu}^\mathrm{GL}),
\qquad n=1,\dots,N,\; \mu=2,\dots,k-1\,.
\]
If $u$ is smooth enough, then the standard Hermite interpolant 
$\It u$ provides the error estimates
\begin{align}
\label{Eq:IntOpI_1}
\| \dt u - \dt \It u\|_{C^0(\overline{I}_n;B)} & \lesssim 
\tau_n^{k+1} \|u\|_{C^{k+2}(\overline{I}_n;B)}\,,\\
\nonumber
\| \dt^2 u - \dt^2 \It u\|_{C^0(\overline{I}_n;B)} & \lesssim 
\tau_n^{k} \|u\|_{C^{k+2}(\overline{I}_n;B)}
\end{align}
on each interval $I_n$.
For a function $u\in C^1(\overline{I};B)$, we construct a local interpolant
$R_n^k u\in \P_k(I_n;B)$ by
\begin{align}
R_n^k u(t_{n-1}^{+}) & = \It u(t_{n-1}^{+})
\intertext{and}
\label{Eq:IntOpR_1}
\dt R_n^k u(t^{\mathrm{GL}}_{n,\mu})
& = \dt \It u(t^{\mathrm{GL}}_{n,\mu}) \,,\quad \mu=1,\dots, k\,,
\end{align}
on each time subinterval $I_n$ and a global interpolant $\Rt u\in
Y_{\tau}^k(B)$ by
\[
(\Rt u)|_{I_n} := R_n^k (u|_{I_n}), \qquad n=1,\dots,N\,.
\]
Finally, we put $\Rt u(0):=u(0)$.

In the following we summarize some basic properties of the operator $\Rt$;
cf.~\cite{BKRS18,ES16} for their proofs.

\begin{lem}
\label{Lem:PropRtau}
Let $u\in C^1(\overline{I};B)$ where $B\subset H$. Then, the function $\Rt u$ is
continuously differentiable in time on $\overline{I}$ with
$\Rt u(t_n) = u(t_n)$ and $\dt R_\tau^{k} u(t_n) = \dt u(t_n)$ for all
$n=0,\dots,N$.

\end{lem}
\begin{lem}
\label{Lem:AppPropR}
For all $n=1,\ldots, N$ and all $u\in C^{k+1}(\overline{I}_n;B)$, there holds that 
\begin{equation}
\label{Eq:AppPropR}
\|u - \Rt u \|_{C^0(\overline{I}_n;B)} \lesssim 
\tau_n^{k+1}\|u\|_{C^{k+1}(\overline{I}_n;B)}\,.
\end{equation}
Moreover, the estimate $\|\Rt u\|_{C^0(\overline{I}_n;B)}\lesssim 
\|u\|_{C^0(\overline{I}_n;B)} + \tau_n \|u\|_{C^1(\overline{I}_n;B)}$ is
satisfied for all $u\in C^1(\overline{I}_n;B)$.
\end{lem}

Lemma \ref{Lem:AppPropR} implies the following result. 

\begin{cor}
\label{Cor:DtRh}
For all $n=1,\ldots, N$ and all $u\in C^{k+1}(\overline{I}_n;B)$, there holds that 
\begin{equation*}
\label{Eq:AppPropDtR}
\|\dt u - \dt \Rt u \|_{C^0(\overline{I}_n;B)} \lesssim 
\tau_n^{k}\|u\|_{C^{k+1}(\overline{I}_n;B)}\,.
\end{equation*}
Moreover, the estimate $\|\dt \Rt u\|_{C^0(\overline{I}_n;B)}\lesssim 
\|u\|_{C^1(\overline{I}_n;B)}$ is satisfied for all
$u\in C^1(\overline{I}_n;B)$.
\end{cor}

For the operator $\IH$ defined in~\eqref{Eq:DefHIntOp} we recall the following 
approximation properties. They directly follow from the standard error estimates 
for Hermite interpolation.

\begin{lem}
\label{Lem:Prop_I_H}
The operator $\IH:C^1(\overline{I};H)\to X^k_\tau(H)$ provides
\begin{align}
\nonumber
\| u - \IH u \|_{C^0(\overline{I}_n;B)} & \lesssim 
\tau_n^{k+1}\|u\|_{C^{k+1}(\overline{I}_n;B)}\,,\\
\label{Eq:AppPropIH_2}
\| \dt u - \dt \IH u \|_{C^0(\overline{I}_n;B)} & 
\lesssim 
\tau_n^{k}\|u\|_{C^{k+1}(\overline{I}_n;B)}\,,
\end{align}
for all $n=1,\ldots, N$ and all $u\in C^{k+1}(\overline{I}_n;B)$.
\end{lem}

Finally, we present a norm bound that will be used later in our
analysis.

\begin{lem}
\label{Lem:HDIR}
For any $u\in \P_k (I_n; H)$ the norm inequality
\begin{equation*}
\label{Lem:HDIR0}
\int_{I_n} \| u \|^2 \ud t \lesssim \tau_n \| u(t_{n-1})\|^2 + \tau_n^2 \int_{I_n} \| 
\dt u \|^2 \ud t
\end{equation*}
holds.
\end{lem}

\section{Error estimates}
\label{Sec:ErrAna}

The overall goal of this work is to prove error estimates for the error
\begin{equation}
\label{Eq:DefTE}
E(t) := U(t) - U_{\tau,h} (t)\,,
\end{equation}
where the Galerkin--collocation approximation $U_{\tau,h}$ is the solution of 
Problem~\ref{Prob:DisLocalcGPC}. We will use in the sequel the componentwise 
representation 
$E(t) = \big(e^0(t),e^1(t)\big)$. We observe that $E$ is continuously differentiable in 
time on $\overline{I}$ if we assume for our analysis that for the exact 
solution $U=\big(u^0,u^1\big) \in \big(C^1(\overline{I};V)\big)^2$ is satisfied.

For each time interval $I_n$, $n=1,\ldots,N$, we define the bilinear form
\begin{equation*}
\label{Def:Bntilde}
\BGL(W,V) := \QGL\big(\llangle 
\dt W, V \rrangle \big)
+ \QGL\big(\llangle \mathcal{A}_h W, V \rrangle \big)
\end{equation*}
where $W$ and $V$ have to satisfy some smoothness conditions to ensure that $\BGL$ is
well-defined.

Our analysis will follow the main lines given in~\cite{BKRS18} since the solution
$U_{\tau,h}$ in this paper is related to $L_{\tau} U_{\tau,h}$ there with the difference
that our polynomial order $k$ is related to $k+1$ in~\cite{BKRS18}. This relation is
motivated by the fact that the solution of the numerically integrated $\cGPone$--$\cGr$,
given in Problem~\ref{Prob:DisLocalcGPC} could also be interpreted as the post-processed
solution of a numerically integrated $\mathrm{cGP(\mathit{k-1})}$--$\cGr$ scheme, given in
Problem~\ref{Prob:SemiDisLoc}, with a modified right-hand side in that $F$ is replaced by
$\IH F$. In order to keep this work self-contained, we will cite the results 
used from~\cite{BKRS18} and will focus on the new aspects in the error analysis.

\subsection{Error estimates for $\boldsymbol{\dt U_{\tau,h}}$}
\label{Sec:ErrDtU}

We start with proving an $L^\infty(L^2)$-norm estimate for the time derivative $\dt E(t)$ of
the error as an auxiliary result. This represents an essential argument in our proof and
is specific to the hyperbolic character of~\eqref{Eq:EP1}. Based on the $L^\infty(L^2)$-bound
for $\dt E(t)$ an estimate for $E(t)$ will be proved in Subsec.~\ref{Sec:L2ErrEst}.

In order to bound $\dt E(t)$, we derive a variational problem that is satisfied by
$\dt U_{\tau,h}$.

\begin{thm}
\label{Thm:EqPtE}
Let $U_{\tau,h}\in \big(X^k_{\tau,h}(V_h)\big)^2$ be the solution 
of Problem~\ref{Prob:DisLocalcGPC}. Then, its time derivative
$\dt U_{\tau,h}\in \big(X^{k-1}_{\tau,h}(V_h)\big)^2$ satisfies for all $n=1,\ldots,N$
the equation
\begin{equation}
\label{Eq:EqPtE0}
\BGL(\dt U_{\tau,h},V_{\tau,h}) = 
\QGL(\llangle \dt \IH F,V_{\tau,h}\rrangle) 
= \int_{I_n} \llangle \dt \IH F,V_{\tau,h}\rrangle \ud t
\end{equation}
for all $V_{\tau,h}\in\big(\P_{k-2}(I_n;V_h)\big)^2$.
\end{thm}
 
\begin{proof}
Recalling that $\dt U_{\tau,h}\in \big(\P_{k-1}(I_n;V_h)\big)^2$, we get by 
the exactness of the Gauss--Lobatto formula \eqref{Eq:GF+GLF} for all polynomials in 
$\P_{2k-3}(I_n;\R)$ along with integration by parts that
\begin{equation}
\label{Eq:EqPtE1_1}
\begin{aligned}
\BGL(\dt U_{\tau,h},V_{\tau,h}) 
&= \QGL(\llangle \underbrace{\dt^2 U_{\tau,h}
+ \mathcal{A}_h \dt U_{\tau,h}}_{\in \mathbb{(}P_{k-1}(I_n;V_h))^2},V_{\tau,h} \rrangle)
= \int_{I_n} \llangle \dt (\dt U_{\tau,h} + \mathcal{A}_h 
U_{\tau,h}),V_{\tau,h}\rrangle \ud t\\
& = -\int_{I_n} \llangle \dt U_{\tau,h} + \mathcal{A}_h 
U_{\tau,h},\dt V_{\tau,h}\rrangle \ud t + \llangle \dt U_{\tau,h} + 
\mathcal{A}_h U_{\tau,h}, V_{\tau,h}\rrangle\Big|_{t_{n-1}^{+}}^{t_n^{-}}
\end{aligned}
\end{equation}
for $V_{\tau,h}\in\big(\P_{k-2}(I_n,V_h)\big)^2$.
Using the exactness of the Hermite quadrature formula $\QH$ for polynomials in
$\P_{2k-3}(I_n;\R)$ and~\eqref{Eq:DisLocalcGPC_4}, we 
conclude from~\eqref{Eq:EqPtE1_1} that 
\begin{equation}
\label{Eq:EqPtE1_2}
\begin{aligned}
\BGL(\dt U_{\tau,h},V_{\tau,h}) 
& = - \QH(\llangle F, \dt V_{\tau,h}\rrangle) + \llangle 
\dt U_{\tau,h} + \mathcal{A}_h U_{\tau,h}, V_{\tau,h}\rrangle
\Big|_{t_{n-1}^{+}}^{t_n^{-}} \\[2ex]
& = - \int_{I_n} \llangle \IH F, \dt V_{\tau,h}\rrangle\ud 
t + \llangle 
\dt U_{\tau,h} + \mathcal{A}_h U_{\tau,h}, V_{\tau,h}\rrangle
\Big|_{t_{n-1}^{+}}^{t_n^{-}} \\[2ex]
& = \int_{I_n} \llangle \dt \IH F, V_{\tau,h}\rrangle\ud t 
- \llangle \IH F, V_{\tau,h} \rrangle
\Big|_{t_{n-1}^{+}}^{t_n^{-}} + 
\llangle \dt U_{\tau,h} + \mathcal{A}_h U_{\tau,h}, V_{\tau,h}\rrangle
\Big|_{t_{n-1}^{+}}^{t_n^{-}}\,.
\end{aligned}
\end{equation}
From~\eqref{Eq:DisLocalcGPC_2} and~\eqref{Eq:DisLocalcGPC_3} along with the interpolation 
properties of $\IH$, it follows that 
\begin{equation}
\label{Eq:EqPtE1_3}
\dt U_{\tau,h}(t_{\ast}) + \mathcal{A}_h U_{\tau,h}(t_{\ast}) = 
\mathcal{P}_h \IH F(t_{\ast})
\end{equation}
for $t_{\ast}\in \{t_{n-1}^{+},t_n^{-}\}$.
Combining~\eqref{Eq:EqPtE1_2} with~\eqref{Eq:EqPtE1_3} 
shows that 
\begin{equation*}
\BGL(\dt U_{\tau,h},V_{\tau,h}) = \int_{I_n} \llangle \dt \IH F, V_{\tau,h}\rrangle\ud t 
\end{equation*}
for all $V_{\tau,h}\in\big(\P_{k-2}(I_n,V_h)\big)^2$.
Recalling that 
$\dt \IH F|_{I_n}\in \big(\P_{k-1}(I_n;V_h)\big)^2$ and the exactness of the Gauss--Lobatto
quadrature for functions of $\P_{2k-3}(I_n;\R)$, this proves the assertion of the 
theorem. 
\end{proof} 
 
\begin{rem}
\label{Rem:VPDtU}
If the solution $u$ of~\eqref{Eq:IBVP} is sufficiently regular, the time derivative
$\dt U = (\dt u, \dt^2 u)$ solves the evolution problem 
\begin{equation}
\label{Eq:EvolPatU}
\dt (\dt U) + \mathcal{A} (\dt U) = \dt F \quad 
\mathrm{in} \;\; (0,T)\,, \qquad \dt U (0) = - \mathcal{A} U(0) + F(0)\,.
\end{equation}
Assumptions on the data such that~\eqref{Eq:EvolPatU} is satisfied can be 
found in, e.g., \cite[p.\ 410, Thm.\ 5]{E10}. 

Rewriting~\eqref{Eq:EqPtE0} as
\begin{equation}
\label{Eq:EqPtE0_1}
\BGL(\dt U_{\tau,h},V_{\tau,h}) = \int_{I_n} \llangle \dt F,V_{\tau,h}\rrangle \ud t
+ \int_{I_n} \llangle \dt \IH F - \dt F,V_{\tau,h}\rrangle \ud t,
\end{equation}
its discrete solution can now be regarded as the $\mathrm{cGP(\mathit{k-1})}$--$\cGr$ approximation of the
evolution problem~\eqref{Eq:EvolPatU} up to the perturbation term
$\int_{I_n} \llangle \dt \IH F - \dt F,V_{\tau,h}\rrangle \ud t$ on the right-hand side.
Further, the collocation
condition~\eqref{Eq:DisLocalcGPC_2} for $n=1$ along with the initial condition that 
$U_{\tau,h}(0)=U_{0,h}$ shows that
$\dt U_{\tau,h}(0) = - \mathcal{A}_h U_{0,h} + \mathcal{P}_h F(0) $ is satisfied.
\end{rem}

We point out that there is a strong analogy between Remark~\ref{Rem:VPDtU} 
and~\cite[Remark\ 5.3]{BKRS18}. The main difference of the two statements comes through  
the different perturbation terms. However,
having in mind the relation of the polynomial orders, both perturbation terms are of the
same approximation order. Hence, we can directly follow the further arguments 
used in~\cite{BKRS18}. Especially, some assumptions about the discrete initial value
$\partial_t U_{\tau,h}(0)$ with respect to the continuous initial value $\partial_t U(0)$
have to be fulfilled.

\begin{lem}
\label{Lem:IniVal}
Let $U_{0,h}:=(R_h u_0, R_h u_1)$. Then there holds that
\begin{equation*}
\label{Eq:IniVal0}
\dt U_{\tau,h} (0) = \begin{pmatrix} R_h & 0 \\ 0 & P_h 
\end{pmatrix}\dt U(0) \,.
\end{equation*}
\end{lem}

We refer to~\cite[Lemma\ 5.4]{BKRS18} for the proof of Lemma~\ref{Lem:IniVal} taking into
consideration that $U_{\tau,h}$ here is associated to $L_{\tau} U_{\tau,h}$
in~\cite{BKRS18}. Also note that the analog of~\cite[Assumption\ 3.6]{BKRS18} is obviously
satisfied by $U_{\tau,h}$ due to~\eqref{Eq:DisLocalcGPC_1} and~\eqref{Eq:DisLocalcGPC_2} for
$n=1$.

Finally, before proving the error estimate for $\partial_t U_{\tau,h}$, we want to cite
Theorem\ 5.5 of~\cite{BKRS18} that is a slightly generalized result of the analysis
in~\cite{KM04} for the $\mathrm{cGP(\mathit{k})}$--$\cGr$ approximation of the wave
equation.

\begin{thm}
\label{Thm:Makridakis}
Let $\ur$ denote the solution of~\eqref{Eq:IBVP} with data 
$\fr$, $\ur_0$, $\ur_1$ instead of $f$, $u_0$, $u_1$. Suppose $\ell\in \mathbb{N}$
and let $\frtau$ be an approximation of $\fr$ such that
\begin{equation}
\label{Eq:ftau_approx}
\| \fr - \frtau \|_{C(\overline{I}_n;H)} \le C_{\hat{f}}\, \tau^{\ell+1}_n \,, 
\qquad n=1, \dots, N,
\end{equation}
where the constant $C_{\hat{f}}$ depends on $\fr$ but is independent of $n$,
$N$, and $\tau_n$. Furthermore, let
$\Urth=\big(\urth^0,\urth^1\big)\in \big(X_{\tau}^\ell(V_h)\big)^2$
be the solution of the local (on $I_n$) perturbed $\mathrm{cGP(\mathit{\ell})}$--$\cGr$
problem
\begin{equation}
\label{Eq:local_probl_hat}
\int_{I_n} \Big(\llangle \dt \Urth , V_{\tau,h} \rrangle 
+ \llangle \mathcal{A}_h \Urth , V_{\tau,h} \rrangle \Big) \ud t = 
\int_{I_n} \llangle \Frtau, V_{\tau,h} \rrangle \ud t
\end{equation}
for all test functions 
$V_{\tau,h} =\big(v_{\tau,h}^0,v_{\tau,h}^1\big) \in \big(\P_{\ell-1}(I_n;V_h)\big)^2$
with $\Frtau :=\big(0, \frtau\big)$
and the initial value $\Urth(t_{n-1}^+) = \Urth(t_{n-1}^-)$ for $n>1$
and $\Urth(t_0)=\widehat{U}_{0,h}:=\big(R_h \ur_0, P_h \ur_1\big)$. For a sufficiently
smooth exact solution $\ur$, the estimates
\begin{align}
\label{Eq:EstErrTD0M}
\| \ur(t) - \urth^0(t) \| + \| \dt\ur(t) - \urth^1(t) \| 
& \lesssim 
\tau^{l+1} \, C_t(\ur) + h^{r+1} \, C_x (\ur) \,,\\
\label{Eq:EstErrTD1M}
\| \nabla\left( \ur(t) - \urth^0(t) \right) \| 
& \lesssim 
\tau^{l+1} \, C_t(\ur) + h^{r} \, C_x (\ur) \,,
\end{align}
hold for all $t\in \overline{I}$ where $C_t(\ur)$ and $C_x(\ur)$ are quantities depending
on various temporal and spatial derivatives of $\ur$.
\end{thm}

We conclude from Theorem~\ref{Thm:Makridakis} the following error estimates.

\begin{thm}
\label{Thm:EstErrTD}
Let $U_{0,h}:=(R_h u_0, R_h u_1)$ and assume that the exact solution
$U = (u^0,u^1) :=(u,\dt u)$ is sufficiently smooth. Then the error estimates
\begin{alignat}{2}
\label{Eq:EstErrTD0}
\| \dt U (t) - \dt U_{\tau,h}(t) \|
& \lesssim  \tau^{k} \, C_t(\dt u) + h^{r+1} \, C_x (\dt u)
& &\lesssim \tau^{k} + h^{r+1} \,,\\
\label{Eq:EstErrTD1}
\|\nabla\left( \dt u^0(t) - \dt u^0_{\tau,h}(t) \right) \| 
& \lesssim \tau^{k} \, C_t(\dt u) + h^{r} \, C_x (\dt 
u) && \lesssim \tau^{k} + h^{r}\,,
\end{alignat}
hold for all $t\in \overline{I}$ where $C_t(\dt u)$ and $C_x(\dt u)$ are quantities
depending on various temporal and spatial derivatives of $\dt u$.
\end{thm}

\begin{proof}
To prove \eqref{Eq:EstErrTD0} and \eqref{Eq:EstErrTD1}, we apply Theorem 
\ref{Thm:Makridakis}. Since the solution $u$ is sufficiently smooth, the function 
$\ur:=\dt u$ is the solution of the wave equation \eqref{Eq:IBVP} with the right-hand 
side $\fr:=\dt f$ and the initial conditions $\ur(0)=\ur_0 := u_1$ and $\dt\ur(0)=\ur_1 
:= f(0) - A u_0$. Let us define the modified right-hand side $\frtau := \dt 
\IH f$
and $\Frtau :=(0,\frtau)$. Then, the discrete function $\Urth := \dt U_{\tau,h}\in 
\big(X_{\tau}^{k-1}(V_h)\big)^2$ satisfies all the conditions required for the discrete
solution $\Urth$ in Theorem~\ref{Thm:Makridakis} with $\ell=k-1$. In fact, by the
construction of the discrete solution $U_{\tau,h}$ in Problem~\ref{Prob:DisLocalcGPC}, 
the continuity of $\dt U_{\tau,h}$ in the discrete points $t_{n}$, $n=0,\ldots, N$, is
ensured by the conditions~\eqref{Eq:DisLocalcGPC_1}--\eqref{Eq:DisLocalcGPC_3}. Therefore,
it holds that $\Urth\in\big(\P_{k-1}(I_n;V_h)\big)^2$ and that 
$\Urth(t_{n-1}^+) = \Urth(t_{n-1}^-) $. Moreover, from 
$U_{0,h}:=\big(R_h u_0, R_h u_1\big)$ and Lemma~\ref{Lem:IniVal}, we get that 
$\widehat{U}_{0,h} = \Urth(0) = \dt U_{\tau,h}(0) = \big(R_h\ur_0, P_h \ur_1\big)$.
Theorem \ref{Thm:EqPtE} implies for all $n=1, \ldots, N$
and all $V_{\tau,h}\in \big(\P_{k-2}(I_n;V_h)\big)^2$ that
\[
\BGL(\Urth, V_{\tau,h}) = \QGL\big(\llangle \dt\Urth,V_{\tau,h}\rrangle + 
\llangle{\mathcal{A}}_h\Urth,V_{\tau,h}\rrangle\big)
= \QGL\big(\llangle \Frtau , V_{\tau,h}\rrangle\big)\,. 
\]
Each quadrature formula in the previous equation is exact since all integrands are 
polynomials in $t$ with degree not greater than $2k-3$ such that the variational equation 
\eqref{Eq:local_probl_hat} of Theorem \ref{Thm:Makridakis} is satisfied. Thus, we have 
shown that $\Urth$ is the discrete solution of Theorem \ref{Thm:Makridakis} for the above 
defined data. To verify the approximation property for $\frtau$, we use the definition of 
$\fr$ and $\frtau$, apply the estimate \eqref{Eq:AppPropIH_2}, and obtain 
\eqref{Eq:ftau_approx} with a constant $C_{\hat{f}}=C \|\dt^{k+1} f\|_{C(\overline{I};H)}$.
Then, we use Theorem \ref{Thm:Makridakis} with $\ell=k-1$. Recalling the 
representation by components, $\dt U = \big(\dt u^0,\dt u^1 \big) = \big(\ur, \dt\ur\big)$
and $\Urth=\big(\urth^0,\urth^1\big) =\big(\dt u^0_{\tau,h},\dt u^1_{\tau,h}\big)$, 
we directly get assertion~\eqref{Eq:EstErrTD0} from~\eqref{Eq:EstErrTD0M}
and assertion~\eqref{Eq:EstErrTD1} from~\eqref{Eq:EstErrTD1M}.
\end{proof}

\subsection{Error estimates for $\boldsymbol{U_{\tau,h}}$}
\label{Sec:L2ErrEst}

This section is devoted to the desired norm estimates for the error 
$E(t) := U(t) - U_{\tau,h}(t)$ where $U_{\tau,h}$ is the solution of
Problem~\ref{Prob:DisLocalcGPC}. For our error 
analysis we consider the decomposition
\begin{equation}
\label{ErrDecomp}
E(t) = \Theta(t) + E_{\tau,h}(t)
\quad\text{with}\quad
\Theta(t) := U(t) - \mathcal{R}_h R_\tau^{k} U(t)\;\text{and}\;
E_{\tau,h} := \mathcal{R}_h R_\tau^{k} U(t)- U_{\tau,h}
\end{equation}
for all $t\in \overline{I}$ and define the components $E_{\tau,h}(t) 
=\big(e_{\tau,h}^{0}(t), e_{\tau,h}^{1}(t)\big)$. We observe that both $\Theta$ and 
$E_{\tau,h}$ are continuously differentiable in time on $\overline{I}$ if the exact
solution $U$ is sufficiently smooth. We refer to $\Theta$ as interpolation error.
We note that both $\Theta$ and $E_{\tau,h}$ are smooth enough to be
used as arguments in the bilinear form $\BGL$.

The following estimates of $\Theta$ in~\eqref{ErrDecomp} can be found
in~\cite[Lemma 5.7]{BKRS18}. They rely on the properties of $\mathcal{R}_h$ and $\Rt$. 

\begin{lem}[Estimation of the interpolation error]
\label{Lem:IntpolErr}
Let $m\in\{0,1\}$. Then, the error estimates
\begin{alignat}{2}
\label{Eq:IntpolErr_2}
\| \Theta (t) \|_m & \lesssim h^{r+1-m} + \tau_n^{k+1},
& \quad t & \in \overline{I}_n\,,\\[1ex]
\label{Eq:IntpolErr_3}
\| \dt\Theta (t) \|_m & \lesssim h^{r+1-m} + \tau_n^{k},
& \quad t & \in \overline{I}_n\,,
\end{alignat}
hold for all $n=1,\ldots, N$ where $\|\cdot\|_0 := \|\cdot\|$.
\end{lem}

Next, we address the discrete error $E_{\tau,h}$ of the decomposition~\eqref{ErrDecomp} 
between the interpolation $\mathcal{R}_h R_\tau^{k} U$ and the fully discrete solution 
$U_{\tau,h}$. We start with some auxiliary results.

\begin{lem}[Consistency error]
\label{Lem:Consist}
Assume that $U\in C^1(\overline{I};V) \times C^1(\overline{I};H)$. Then, for all 
$n=1,\ldots, N$ the identity 
\begin{equation*}
\BGL (E, V_{\tau,h}) = \QGL \big( \llangle \IGL F- \IH F, V_{\tau,h} \rrangle \big)
= \QGL \big( \llangle F- \IH F, V_{\tau,h} \rrangle\big)
\end{equation*}
is satisfied for all $V_{\tau,h}\in \big(Y_{\tau,h}^{k-2}(V_h)\big)^2$.
\end{lem}

\begin{proof}
We recall from Lemma~\ref{Lem:DSVP} that for all $n=1,\ldots, N$ the identity
\begin{equation}
\label{Eq:Cons_1}
\BGL \big(U_{\tau,h},V_{\tau,h}\big) 
= \QGL \big(\llangle \IH F,V_{\tau,h}\rrangle\big)
\end{equation}
holds for all $V_{\tau,h}\in \big(\P_{k-2}(I_n;V_h)\big)^2$. We have
under sufficient smoothness assumptions on the exact solution that 
\begin{equation}
\label{Eq:Cons_2}
\dt U(t^{\mathrm{GL}}_{n,\mu}) + \mathcal{A}U (t^{\mathrm{GL}}_{n,\mu})= F 
(t^{\mathrm{GL}}_{n,\mu}), \qquad \mu=1,\ldots,k\,.
\end{equation}
By the consistency~\eqref{Eq:ConsistA_h} of $\mathcal{A}_h$, the identity~\eqref{Eq:Cons_2}
implies
\begin{equation}
\label{Eq:Cons_3}
\begin{aligned}
\BGL\big(U,V_{\tau,h}\big) & = \QGL\big(\llangle \dt U + \mathcal{A}_h U,
V_{\tau,h}\rrangle\big)\\
& = \QGL\big(\llangle \dt U + \mathcal{A} U, V_{\tau,h}\rrangle\big) 
= \QGL\big(\llangle \IGL F, V_{\tau,h}\rrangle\big) \,.
\end{aligned}
\end{equation}
Combining~\eqref{Eq:Cons_1} with~\eqref{Eq:Cons_3} and recalling that $E = U 
- U_{\tau,h}$ prove the assertion.
\end{proof}

The following lemma is slightly more general than~\cite[Lem.~5.9]{BKRS18} where the
proof can be found.
\begin{lem}
\label{Lem:TdTE}
Let $p\in\P_k(I_n)$ be an arbitrary polynomial of degree less than or equal to $k$.
Then, the relation
\begin{equation*}
\label{Eq:TdTE:0}
\dt p(t_{n,\mu}^{\mathrm{G}}) = \dt \IGL p(t_{n,\mu}^{\mathrm{G}})
\end{equation*}
holds for all Gauss points $t^G_{n,\mu}\in I_n$, $\mu=1,\ldots,k-1$.
\end{lem}

Exploiting the correspondence of $U_{\tau,h}$ in this paper to $L_{\tau}U_{\tau,h}$
in~\cite{BKRS18} and keeping in mind that $k$ here is related to $k+1$ there, we can recall
from~\cite{BKRS18} the results of stability (cf.\ \cite[Lemma 5.10]{BKRS18}) and 
boundedness (cf.\ \cite[Lemma 5.11]{BKRS18}).

\begin{lem}[Stability]
\label{Lem:Stab_2} 
We have
\begin{multline}
\label{Eq:StabEst_21}
\BGL \big( (e_{\tau,h}^{0},e_{\tau,h}^{1}), 
(\Pin A_h \IGL e_{\tau,h}^{0}, \Pin \IGL e_{\tau,h}^{1} )\big)\\
= \frac{1}{2}\left( \|\nabla e_{\tau,h}^{0}(t_n)\|^2 
-\|\nabla e_{\tau,h}^{0}(t_{n-1})\|^2 + \| e_{\tau,h}^{1}(t_n)\|^2 - \| 
e_{\tau,h}^{1}(t_{n-1})\|^2 \right)
\end{multline}
for all $n=1,\ldots, N$.
\end{lem}

\begin{lem}[Boundedness]
\label{Lem:Bound}
Let $V_{\tau,h} = \big(\Pin A_h \IGL e_{\tau,h}^{0}, 
\Pin \IGL e_{\tau,h}^{1} \big)$. Then, the bound
\begin{equation*}
\left| \BGL \big(\Theta, V_{\tau,h})\right|
\lesssim \tau_n^{1/2} \big(\tau_n^{k+1} + h^{r+1}\big)
\left\{ \tau_n \| E_{\tau,h}(t_{n-1})\|^2 + \tau_n^2 
\QG\big(\| \dt E_{\tau,h}\|^2\big) \right\}^{1/2}
\end{equation*}
holds for all $n=1,\ldots, N$.
\end{lem}

We proceed with estimating the consistency error given in Lemma~\ref{Lem:Consist}.

\begin{lem}[Estimates on right-hand side term]
\label{Lem:EstRhs}
Let $V_{\tau,h} = \big(v_{\tau,h}^{0},v_{\tau,h}^{1} \big)
=\big(\Pin A_h \IGL e_{\tau,h}^{0},\Pin$ $\IGL e_{\tau,h}^{1} \big)$.
Then, the estimate
\begin{equation*}
\QGL\big(\llangle (0, f - \IH f), (v_{\tau,h}^{0},v_{\tau,h}^{1}) \rrangle\big)
\lesssim \tau_n^{1/2} \tau_n^{k+1}
\left\{ \tau_n \| E_{\tau,h}(t_{n-1})\|^2 + \tau_n^2 
\QG \big(\| \dt E_{\tau,h}\|^2\big) \right\}^{1/2}
\end{equation*}
holds for all $n=1,\ldots, N$.
\end{lem}

\begin{proof}
Using the Cauchy--Schwarz inequality along with Lemma~\ref{Lem:Prop_I_H}, we get that 
\begin{align*}
\QGL\big(\llangle (0, f - \IH f), (v_{\tau,h}^{0},v_{\tau,h}^{1}) \rrangle\big)
& = \QGL\big(\langle f- \IH f, \Pin \IGL e_{\tau,h}^{1} \rangle\big)\\[1ex]
& \leq \big(\QGL\big(\| f-\IH f \|^2\big)\big)^{1/2}
\left(\QGL\big(\| \Pin \IGL e_{\tau,h}^{1} \|^2\right)\big)^{1/2}\\[1ex]
& \lesssim \tau_n^{1/2}\tau_n^{k+1}
\left(\QGL\big(\| \Pin \IGL e_{\tau,h}^{1} \|^2\right)\big)^{1/2}\,.
\end{align*}
Using the exactness of $\QGL$ for polynomials up to degree $2k-3$, the
stability of the $L^2$-projection $\Pin$, the norm bound from
Lemma~\ref{Lem:HDIR}, and Lemma~\ref{Lem:TdTE}, we finally conclude that
\begin{align*}
\QGL\big(\|\Pin \IGL e_{\tau,h}^1\|^2 \big)
= \int_{I_n} \|\Pin \IGL e_{\tau,h}^1\|^2 \ud t
& \le \int_{I_n} \|\IGL e_{\tau,h}^1\|^2 \ud t\\
& \lesssim \tau_n \|\IGL e_{\tau,h}^1(t_{n-1})\|^2
  + \tau_n^2 \int_{I_n} \|\dt \IGL e_{\tau,h}^1\|^2\ud t\\[1ex]
& = \tau_n \|e_{\tau,h}^1(t_{n-1})\|^2 + \tau_n^2 \QG\big( \|\dt \IGL e_{\tau,h}^1\|^2 
\big)\\[1ex]
& = \tau_n \|e_{\tau,h}^1(t_{n-1})\|^2 + \tau_n^2 \QG\big( \|\dt e_{\tau,h}^1\|^2 \big)\,.
\end{align*}
Combining both estimates, the assertion of the lemma follows directly.
\end{proof}

\begin{lem}[Estimates on $E_{\tau,h}$]
\label{Lem:EstTEth}
Let $U_{0,h}:=\big(R_h u_0, R_h u_1\big)$. Then, the estimate
\begin{equation}
\label{Eq:TE_01}
\| e_{\tau,h}^{0}(t_n)\|_1^2 + \| e_{\tau,h}^{1}(t_n)\|^2 
\lesssim \big(\tau^{k+1}+h^{r+1}\big)^2
\end{equation}
is satisfied for all $n=1,\ldots,N$. Moreover, we have that 
\begin{align}
\label{Eq:TE_012}
\|\nabla e_{\tau,h}^{\; 0}(t)\| & \lesssim \tau^{k+1}+h^{r}\,,\\[1ex]
\label{Eq:TE_013}
\| e_{\tau,h}^{\; 0}(t)\| + \| e^{1}_{\tau,h}(t)\|
& \lesssim \tau^{k+1}+h^{r+1} 
\end{align}
for all $t\in\overline{I}$.
\end{lem}

\begin{proof}
We conclude from Lemma~\ref{Lem:Consist} that
\begin{equation*}
\BGL \big(E_{\tau,h}, V_{\tau,h}\big) = -\BGL \big(\Theta, V_{\tau,h}\big)
+ \QGL \big( \llangle F- \IH F, V_{\tau,h} \rrangle\big)
\end{equation*}
is satisfied for all $V_{\tau,h}\in \big(Y_{\tau,h}^{k-2}(V_h)\big)^2$.
Choosing here $V_{\tau,h} = \big(\Pin A_h \IGL e_{\tau,h}^{0}, \Pin \IGL e_{\tau,h}^{1}\big)$
and using Lemma~\ref{Lem:Bound} and Lemma~\ref{Lem:EstRhs} yield that 
\begin{multline}
\label{Eq:TE_03}
\BGL \big( (e_{\tau,h}^{0}, e_{\tau,h}^{1}), (\Pin A_h \IGL e_{\tau,h}^{0}, 
\Pin \IGL e_{\tau,h}^{1} )\big)\\
\begin{aligned}
& \qquad = - \BGL \big( (\theta^0, \theta^1), (\Pin A_h \IGL e_{\tau,h}^{0},
\Pin \IGL e_{\tau,h}^{1} )\big)\\[1ex] 
& \qquad \quad + \QGL \big( (0, f - \IH f), (\Pin A_h \IGL e_{\tau,h}^{0},
\Pin \IGL e_{\tau,h}^{1} ) \big)\\[1ex]
& \qquad \lesssim 
\tau_n^{1/2}\big(\tau_n^{k+1}+h^{r+1}\big)
\left\{ \tau_n \| E_{\tau,h}(t_{n-1})\|^2 + 
\tau_n^2 Q_{n}^{\mathrm G}(\| \dt E_{\tau,h}\|^2) \right\}^{1/2} \,.
\end{aligned}
\end{multline}
Since the upper bound in~\eqref{Eq:TE_03} coincides with that in Eq.~(5.46) of~\cite{BKRS18}
and our $E_{\tau,h}$ can be identified with $\widetilde{E}_{\tau,h}$ of~\cite{BKRS18}, we present
here just a short summary of the proof of Lemma~5.12 in~\cite{BKRS18}.

Combining the stability property~\eqref{Eq:StabEst_21} of $\BGL$ with~\eqref{Eq:TE_03}, applying the
Cauchy--Schwarz inequality, and telescopic summing lead to
\begin{align}
\|\nabla e_{\tau,h}^{0}(t_n)\|^2 + \| e_{\tau,h}^{1}(t_n)\|^2 
& \lesssim \|\nabla e_{\tau,h}^{0}(t_0)\|^2 + \| e_{\tau,h}^{1}(t_0)\|^2
+ \sum_{s=1}^n \tau_s (\tau_s^{k+1}+h^{r+1})^2 \nonumber\\
&\qquad + \sum_{s=1}^n \tau_s^2 
Q_s^{\mathrm{G}} \big(\| \dt E_{\tau,h}\|^2 \big) + \sum_{s=1}^n \tau_s 
\| E_{\tau,h}(t_{s-1})\|^2.
\label{Eq:TE_05}
\end{align}
Using 
\begin{equation}
\label{Eq:TE_05a}
\| \dt E_{\tau,h}(t)\| \le 
\| \dt U (t) - \dt U_{\tau,h}(t) \| +
\|-\dt\Theta(t)\|
\lesssim \tau^{k} + h^{r+1}\,, \qquad t\in\overline{I}\,,
\end{equation}
together with the estimates~\eqref{Eq:EstErrTD0} and~\eqref{Eq:IntpolErr_3}, we obtain 
that
\begin{equation}
\label{Eq:TE_06}
\begin{aligned}
\|\nabla e_{\tau,h}^{0}(t_n)\|^2 + \| e_{\tau,h}^{1}(t_n)\|^2 
& \lesssim \|\nabla e_{\tau,h}^{0}(t_0)\|^2 + \| e_{\tau,h}^{1}(t_0)\|^2
+ (\tau^{k+1}+h^{r+1})^2 \\
&\qquad
+ {\sum_{s=0}^{n-1} \tau_{s+1}} (\|\nabla 
e_{\tau,h}^{\; 0}(t_{s})\|^2 +\| e^{1}_{\tau,h}(t_{s})\|^2) \,,
\end{aligned}
\end{equation}
where we also used the definition of the Gauss quadrature and the Poincar\'e 
inequality. Applying the discrete Gronwall lemma (cf.\ \cite[p.\ 14]{Q08}) results in
\begin{equation*}
\|\nabla e_{\tau,h}^{0}(t_n)\|^2 + \| e_{\tau,h}^{1}(t_n)\|^2 
\lesssim \|\nabla e_{\tau,h}^{0}(t_0)\|^2 + \| e_{\tau,h}^{1}(t_0)\|^2 + 
(\tau^{k+1}+h^{r+1})^2 \,.
\end{equation*}
Exploiting $e_{\tau,h}^{i}(t_0) = 0$, $i\in\{0,1\}$, which holds due to the
choice $U_{0,h}=\big(R_h u_0, R_h u_1\big)$ of the discrete initial value,
this estimate along with the Poincar\'e inequality proves the
assertion~\eqref{Eq:TE_01}.

To show~\eqref{Eq:TE_012} and~\eqref{Eq:TE_013}, we start for the error
component $ e^i_{\tau,h}\in \mathbb{P}_{k}(I_n,V_h)$, $i\in \{0,1\}$, with 
\begin{equation}
\label{Eq:TE_max_e}
\| e^i_{\tau,h}(t)\|_m \le \| e^i_{\tau,h}(t_{n})\|_m 
+
\tau_n \max_{s\in\overline{I}_n} \|\dt e^i_{\tau,h}(s)\|_m\,, 
\qquad t\in I_n\,,
\end{equation}
that is deduced from the fundamental theorem of calculus. Applying~\eqref{Eq:TE_01} 
and~\eqref{Eq:TE_05a}, we get from~\eqref{Eq:TE_max_e} with $m=0$ that
\[
\| e^i_{\tau,h}(t)\| \lesssim (\tau^{k+1} + h^{r+1})
+ \tau_n (\tau^{k} + h^{r+1}) \lesssim \tau^{k+1} + h^{r+1}\,,
\qquad t\in\overline{I},\,i\in \{0,1\}\,,
\]
which proves \eqref{Eq:TE_013}.

Similarly to~\eqref{Eq:TE_05a}, we get for the $H^1$-norm that 
\begin{equation}
\label{Eq:TE_05b}
\begin{aligned}
\| \dt e^0_{\tau,h}(t)\|_1 & \le 
\| \dt u^0 (t) - \dt u^0_{\tau,h}(t) \|_1 +
\|-\dt\theta^0(t)\|_1
\lesssim \tau^{k} + h^{r}\,, \qquad t\in\overline{I}\,,
\end{aligned}
\end{equation}
where we used~\eqref{Eq:EstErrTD1} along with the Poincar\'e inequality 
and~\eqref{Eq:IntpolErr_3}. 
Applying~\eqref{Eq:TE_01} and~\eqref{Eq:TE_05b}, we get from~\eqref{Eq:TE_max_e} with $m=1$
that
\[
\| e^0_{\tau,h}(t)\|_1 \lesssim (\tau^{k+1} + h^{r+1})
+ \tau_n (\tau^{k+1} + h^{r}) \lesssim \tau^{k+1} + h^{r}\,,\qquad t\in\overline{I}\,,
\]
which proves~\eqref{Eq:TE_012}.
\end{proof}

We are now able to derive our final error estimates for the proposed 
Galerkin--collocation approximation of the solution to~\eqref{Eq:IBVP}.

\begin{thm}[Error estimate for $U_{\tau,h}$]
\label{Thm:OvEst}
Let $U=(u,\dt u)$ be the solution of the problem~\eqref{Eq:IBVP} and let $U_{\tau,h}$ be
the fully discrete solution of Problem~\ref{Prob:DisLocalcGPC} with initial value
$U_{0,h} = (R_h u_0, R_h u_1)$. Then, the error $ E(t) = \big(e^0(t), e^1(t)\big)
= U(t) - U_{\tau,h}(t)$ can be bounded for all $t\in\overline{I}$ by
\begin{align}
\label{Eq:ErrE_01}
\| e^0(t) \| + \| e^1(t) \| & \lesssim \tau^{k+1}+h^{r+1}\,,\\[1ex]
\label{Eq:ErrE_01a}
\| \nabla e^0(t) \| & \lesssim \tau^{k+1}+h^{r}\,.
\end{align}
Moreover, the estimates
\begin{align}
\label{Eq:ErrE_02}
\| e^0\|_{L^2(I;H)} + \| e^1 \|_{L^2(I;H)} & \lesssim \tau^{k+1}+h^{r+1}\,,\\[1ex]
\label{Eq:ErrE_022}
\| \nabla e^0\|_{L^2(I;H)} & \lesssim \tau^{k+1}+h^{r}
\end{align}
hold true.
\end{thm}

\begin{proof}
Recalling the error decomposition 
\begin{equation}
\label{Eq:TotES_0}
 E(t) = U(t) - U_{\tau,h}(t) = \Theta(t) + E_{\tau,h}(t) \,, 
\end{equation}
we conclude assertion~\eqref{Eq:ErrE_01} by applying the triangle inequality along 
with estimate~\eqref{Eq:IntpolErr_2} with $m=0$ and~\eqref{Eq:TE_013} to the terms on the 
right-hand-side of~\eqref{Eq:TotES_0}. 
Similarly we conclude~\eqref{Eq:ErrE_01a} using the estimate~\eqref{Eq:IntpolErr_2} with
$m=1$ and~\eqref{Eq:TE_012}. The assertions~\eqref{Eq:ErrE_02} and~\eqref{Eq:ErrE_022}
follow from the definition of the $L^2(I;H)$-norm together with the
estimates~\eqref{Eq:ErrE_01} and~\eqref{Eq:ErrE_01a}.
\end{proof}

\begin{rem}
We note that the estimates \eqref{Eq:ErrE_01} to \eqref{Eq:ErrE_022} are of optimal 
order in space and time. 
\end{rem}

Similarly to the estimate of $\partial_t E_{\tau,h}$ in Subsec.~\ref{Sec:ErrDtU}, the 
estimation of $E_{\tau,h}$ in Lemma~\ref{Lem:EstTEth} follows the analysis of 
\cite{BKRS18}. The main difference comes through the consistency error of 
Lemma~\ref{Lem:Consist}. However, this  does not cause any difficulties since the 
consistency error is of the same order as the further terms that get involved in the 
error analysis through Lemma~\ref{Lem:Bound} of boundedness.

\section{Energy conservation principle for $\boldsymbol{f\equiv 0}$}
\label{Sec:EngCons}

In this section we address the issue of energy conservation for the considered space-time finite 
element scheme. For vanishing right-hand side term $f\equiv 0$ it is well-known 
that the solution $u$ of the initial-boundary value problem \eqref{Eq:IBVP} satisfies 
the energy conservation 
\begin{equation*}
\| u^1(t) \|^2 + \| \nabla u^0(t) \|^2 = 
\| u_1 \|^2 + \| \nabla u_0 \|^2, \qquad
t\in I. 
\end{equation*}

We will prove that the space-time finite element discretization $U_{\tau,h}$ of
Problem~\ref{Prob:DisLocalcGPC} also satisfies the energy conservation principle at the
discrete time nodes $t_n$. Preserving this fundamental property of the solution
of~\eqref{Eq:IBVP} is an important quality criterion for discretization schemes
of~\eqref{Eq:IBVP}. 

\begin{lem}[Energy conservation for $U_{\tau,h}$]
Suppose that $f\equiv 0$. Let the initial value be given by 
$U_{0,h}=(u_{0,h},u_{1,h})$. Then, the fully discrete solution 
$U_{\tau,h}=(u^0_{\tau,h},u^1_{\tau,h})$ defined by 
Problem~\ref{Prob:DisLocalcGPC} satisfies the energy conservation property
\begin{equation}
\label{Eq:EngConsU}
\| u^1_{\tau,h}(t_n)\|^2 + \| \nabla u^0_{\tau,h}(t_n) \|^2 =
\| u_{1,h} \|^2 + \| \nabla u_{0,h} \|^2
\end{equation}
for all $n=1,\ldots,N$.
\end{lem}

\begin{proof}
Let $f\equiv 0$. We recall that the fully discrete 
solution $U_{\tau,h}=\big(u^0_{\tau,h},u^1_{\tau,h}\big)$ defined by 
Problem~\ref{Prob:DisLocalcGPC} satisfies the variational
equation~\eqref{Lem:DSVP_0}. We choose the test function
$V_{\tau,h} = \big(-\dt \IGL u^1_{\tau,h},
\dt \IGL u^0_{\tau,h}\big)\in \big(\P_{k-2}(I_n;V_h)\big)^2$.
Then, we get from the definitions of $\IGL$ and $\QGL$ that
\begin{align*}
0 & = \QGL\big(\llangle (\dt u^0_{\tau,h}, \dt u^1_{\tau,h} ),
(-\dt \IGL u^1_{\tau,h}, \dt \IGL u^0_{\tau,h} ) \rrangle\big) \\[1ex]
& \quad + \QGL \big(\llangle (- u^1_{\tau,h}, A_h u^0_{\tau,h}),
(-\dt \IGL u^1_{\tau,h}, \dt \IGL u^0_{\tau,h} ) \rrangle\big) \\[2ex]
& = \QGL \big (\llangle (\dt u^0_{\tau,h}, \dt u^1_{\tau,h} ),
(-\dt \IGL u^1_{\tau,h}, \dt \IGL u^0_{\tau,h} ) \rrangle \big) \\[1ex]
& \quad + \QGL \big( \llangle (- \IGL u^1_{\tau,h}, A_h \IGL u^0_{\tau,h}),
(-\dt \IGL u^1_{\tau,h}, \dt \IGL u^0_{\tau,h} ) \rrangle \big).
\end{align*}
Setting
\begin{align*}
T_1 & := \int_{I_n} \llangle (\dt u^0_{\tau,h}, \dt u^1_{\tau,h}),
(-\dt \IGL u^1_{\tau,h}, \dt \IGL u^0_{\tau,h})\rrangle\ud t, \\[1ex]
T_2 & := \int_{I_n} \llangle (- I_\tau^{\mathrm{GL}} u^1_{\tau,h},A_h \IGL u^0_{\tau,h}),
(-\dt \IGL u^1_{\tau,h}, \dt \IGL u^0_{\tau,h} ) \rrangle \ud t,
\end{align*}
the exactness of the $k$-point Gauss--Lobatto quadrature for all polynomials of
maximum degree $2k-3$ gives now
\begin{equation}
\label{Eq:EngConsU_0} 
0 = T_1 + T_2.
\end{equation}
We conclude for $T_1$ by using the exactness of the $(k-1)$-point Gauss quadrature for all
polynomials of maximum degree $2k-3$ and Lemma~\ref{Lem:TdTE} that
\begin{equation}
\label{Eq:EngConsU_1}
\begin{aligned}
T_1 & = \QG\big(\llangle (\dt u^0_{\tau,h}, \dt u^1_{\tau,h}),
(-\dt \IGL u^1_{\tau,h}, \dt \IGL u^0_{\tau,h} ) \rrangle\big)\\[1ex]
& = \QG \big(\llangle (\dt u^0_{\tau,h}, \dt u^1_{\tau,h} ),
(-\dt u^1_{\tau,h}, \dt u^0_{\tau,h} )\rrangle\big) = 0.
\end{aligned}
\end{equation}
Recalling~\eqref{Eq:DefAh}, it follows for $n=1,\ldots,N$ that
\begin{equation}
\label{Eq:EngConsU_2}
\begin{aligned}
T_2 & = \int_{t_{n-1}}^{t_n} \left( \frac{1}{2}
\left\{\text{d}_t \|\IGL u^1_{\tau,h}\|^2 +
\text{d}_t \|\nabla \IGL u^0_{\tau,h}\|^2 \right\}\right)\ud t\\[1ex]
& = \frac{1}{2}\big(\|\IGL u^1_{\tau,h}(t_n)\|^2 - \|\IGL u^1_{\tau,h}(t_{n-1})\|^2
+ \|\nabla \IGL u^0_{\tau,h}(t_n)\|^2 
- \|\nabla \IGL u^0_{\tau,h}(t_{n-1})\|^2 \big).
\end{aligned}
\end{equation}
Now, we combine~\eqref{Eq:EngConsU_0} with~\eqref{Eq:EngConsU_1} and~\eqref{Eq:EngConsU_2}.
We change in the resulting identity the index $n$ to $m$ and sum up from $m=1$ to $n$.
We recall that $\IGL U_{\tau,h}(t_{\ast}) = U_{\tau,h}(t_{\ast})$ for
$t_\ast \in \{t_{n-1},t_n\}$ by definition of $\IGL$ along with
$U_{\tau,h} \in \big(C(\overline{I};V_h)\big)^2$. Hence, assertion~\eqref{Eq:EngConsU}
follows directly.
\end{proof}

\section{C$^{\boldsymbol 2}$-regular Galerkin-collocation approximation and 
its relation to post-processed cGP--C$^{\boldsymbol 1}$}
\label{Sec:PP}

In this section, let $k\ge 5$ be satisfied. Firstly, we propose a family of
Galerkin--collocation time discretization schemes with twice continuously
differentiable in time discrete solutions, that are referred to as
$\cGPtwo$--$\cGr$ schemes. Similarly to the $\cGPone$--$\cGr$ approach of
Problem~\ref{Prob:DisLocalcGPC}, the higher order regularity in time is
ensured by collocation conditions that are imposed in the endpoints
$t_{n-1}$ and $t_n$ of the subinterval $I_n$. This construction principle
can be generalized to discrete solutions of even higher order regularity in
time. For this generalization we also refer to~\cite{BM19,BMW17} where the
Galerkin--collocation approximation of first-order ordinary differential
equations systems is studied in detail. Secondly, we show how the
cGP--C$^2$($k+1$)--cG($r$) approximation can be computed  efficiently in a
simple and computationally cheap post-processing step from the
$\cGPone$--$\cGr$ approach. The post-processing introduced in~\cite{MS11}
and generalized in~\cite{BM19} was recently applied in~\cite{BKRS18} to the
cGP($k$)--cG($r$) family of schemes given in Problem~\ref{Prob:SemiDis}.
There the post-processing is used to lift continuous in time
discrete solutions to continuously differentiable ones. Moreover,
an optimal order error analysis is provided for the post-processed
solution. 

\begin{prob}[Local, numerically integrated, fully discrete problem of $\cGPtwo$--$\cGr$
on $I_n$]
\label{Prob:PP}
\mbox{}\\
Given $U_{\tau,h}(t_{n-1}^-)$ for $n>1$ and $U_{\tau,h}(t_0^-)=U_{0,h}$ for $n=1$,
find $U_{\tau,h}|_{I_n} \in \big(\P_k(I_n;V)\big)^2$ such that
\begin{subequations}
\begin{align}
\label{Eq:PP_1}
U_{\tau,h}(t_{n-1}^+) & = U_{\tau,h}(t_{n-1}^-) \,,\\[1ex]
\label{Eq:PP_2}
\dt U_{\tau,h}(t_{n-1}^+) & = - \mathcal{A}_h U_{\tau,h}(t_{n-1}^+)
+ \mathcal{P}_h F(t_{n-1}^+) \,, \\[1ex]
\label{Eq:PP_2D}
\dt^2 U_{\tau,h}(t_{n-1}^+) & = - \mathcal{A}_h \dt U_{\tau,h}(t_{n-1}^+)
+ \mathcal{P}_h \dt F(t_{n-1}^+) \,, \\[1ex]
\label{Eq:PP_3}
\dt U_{\tau,h}(t_{n}^-) & = - \mathcal{A}_h U_{\tau,h}(t_{n}^-)
+ \mathcal{P}_h F(t_{n}^-) \,, \\[1ex]
\label{Eq:PP_3D}
\dt^2 U_{\tau,h}(t_{n}^-) & = - \mathcal{A}_h \dt U_{\tau,h}(t_{n}^-)
+ \mathcal{P}_h \dt F(t_{n}^-) \,,
\intertext{and}
\label{Eq:PP_4}
Q^\mathrm{H}_{n,k} \Big(\llangle \dt U_{\tau,h} , V_{\tau,h} \rrangle 
& + \llangle \mathcal{A}_h U_{\tau,h} , V_{\tau,h} \rrangle \Big) = 
Q^\mathrm{H}_{n,k} \big(\llangle F,V_{\tau,h}\rrangle\big)
\end{align}
\end{subequations}
for all $V_{\tau,h} \in \big(\P_{k-5} (I_n;V_h)\big)^2$.
\end{prob}
We note that a Hermite-type quadrature formula with $k$ evaluations of function 
values is used in \eqref{Eq:PP_4}. This differs from $\QH$ in \eqref{Eq:GLHF} 
that is used in the $\cGPone$--$\cGr$ family of schemes of 
Problem~\ref{Prob:DisLocalcGPC} and is based on $k-1$ evaluations of function values 
only. In both cases the derivatives of the integrand are evaluated additionally in the 
endpoints of the subinterval $I_n$. Further, the $\cGPtwo$ approach presented here 
differs from that in~\cite{BMW17} by the applied quadrature formula.

\begin{rem}
A careful inspection of the conditions on $U_{\tau,h}$ shows that
\[
\dt U_{\tau,h}(t_{n-1}^+) = \dt U_{\tau,h}(t_{n-1}^-)
\qquad\text{and}\qquad
\dt^2 U_{\tau,h}(t_{n-1}^+) = \dt^2 U_{\tau,h}(t_{n-1}^-)\,,
\]
where the discrete initial conditions are determined using
\[
\dt U_{\tau,h}(0) = - \mathcal{A}_h U_{\tau,h}(0)
+ \mathcal{P}_h F(0), \qquad
\dt^2 U_{\tau,h}(0) = - \mathcal{A}_h \dt U_{\tau,h}(0)
+ \mathcal{P}_h \dt F(0).
\]
Hence, the obtained trajectory in time is twice continuously
differentiable on $\overline{I}$.

Compared to Problem~\ref{Prob:DisLocalcGPC}, the test space of 
the condition~\eqref{Eq:PP_4} is decreased from $\big(\P_{k-3} (I_n;V_h)\big)^2$
to $\big(\P_{k-5} (I_n;V_h)\big)^2$ while the number of collocation conditions
is increased from two to four. For $k=5$ this results in a test space which consists of
piecewise constant functions only and to two additional collocation conditions in both 
endpoints of the time subinterval $I_n$.
\end{rem}

%
%

Finally we address the connection between the ${\mathrm{cGP\text{-}C^1}}$
and ${\mathrm{cGP\text{-}C^2}}$ families of Galerkin--collocation schemes.
\begin{thm}
\label{Thm:PP}
Let $U_{\tau,h}$ denote the solution of the $\cGPone$--$\cGr$ method given 
in Problem~\ref{Prob:DisLocalcGPC}. For $n=1,\dots,N$ we put 
\[
\widetilde{U}_{\tau,h}|_{I_n} := U_{\tau,h}|_{I_n} - K_n \vartheta_n\,,
\]
where $\vartheta_n\in\P_{k+1}(I_n;\mathbb{R})$ is uniquely determined by
\[
I_n^\mathrm{H}\vartheta_n \equiv 0
\qquad\text{and}\qquad
\dt^2 \vartheta_n (t_{n-1}^+) = 1\,.
\]
If the correction coefficient $K_n$ is chosen as
\[
K_n := \begin{cases}
\dt^2 U(t_0^+) - \dt^2 u(t_0)\,, & n=1\,, \\
\dt^2 U(t_{n-1}^+) - \dt^2 \widetilde{U}(t_{n-1}^-)\,, & n>1\,,
\end{cases}
\]
then $\widetilde{U}_{\tau,h}\in \big(X^{k+1}_{\tau}(V_h)\big)^2$ is the solution of the
$\mathrm{cGP\text{-}C^2(\mathit{k}+1)}$--$\cGr$ method given in
Problem~\ref{Prob:PP}.
\end{thm}
The post-processing or lifting operator that is introduced in
Theorem~\ref{Thm:PP} is similar to the lifting operator of~\cite{BKRS18}
that is studied there in the context of the cGP($k$)-cG($r$) approach of
Problem~\ref{Prob:SemiDisLoc}. Both post-processing procedures provide the
correction as a product of a scalar polynomial $\vartheta_n$ and a
coefficient $K_n\in V_h^2$ that are, however, different for the two
procedures. In particular, the lifting in~\cite{BKRS18} is based on the
difference of first derivatives while our post-processing uses the difference
of second order derivatives. We refer to~\cite{BM19} for details on
post-processing techniques for general nonlinear systems of ordinary
differential equations and the proof of the analogue to
Theorem~\ref{Thm:PP}.

\section{Numerical studies}
\label{Sec:NumExp}

In this section we present the results of two numerical experiments for 
the Galerkin--collocation approximation schemes introduced 
in Problem~\ref{Prob:DisLocalcGPC} and Problem~\ref{Prob:PP}, respectively. In 
particular, we aim to illustrate the error estimates given in Theorem~\ref{Thm:OvEst} for 
the $\cGPone$--$\cGr$ Galerkin--collocation approximation 
of Problem~\ref{Prob:DisLocalcGPC}. The 
implementation of the numerical schemes was done in the high-performance 
\texttt{DTM++/awave} frontend solver (cf.\ \cite{K15}) for the \texttt{deal.II} library 
\cite{DealIIReference}. For further details including a presentation of the 
applied algebraic solver and preconditioner we refer to \cite{AB19,K15}.

\subsection{Convergence test for
$\boldsymbol{\mathrm{cGP\text{-}C^1(3)}}$--$\boldsymbol{\mathrm{cG}(3)}$}

%
\begin{table}[htb!]
\centering
\caption{%
Calculated errors $E = \big(e^0, e^1\big)$ with $E(t) = U(t) - U_{\tau,h}(t)$ 
and corresponding experimental orders of convergence (EOC) for the solution $U = (u, 
\dt u)$ of \eqref{Eq:ExSol1} and the Galerkin--collocation approximation 
$U_{\tau,h}\in \left(X_\tau^3 \big(V_h^{(3)}\big)\cap C^1\big(I;V_h^{(3)}\big)\right)^2$ 
of 
Problem \ref{Prob:DisLocalcGPC}.
}
\label{Tab:1}
{
\begin{tabular}{cccccccc}
\toprule
{$\tau$} & {$h$} &
{ $\| e^{0} \|_{L^\infty(L^2)} $ } &
{ $\| e^{1} \|_{L^\infty(L^2)} $ } &
{ $||| E |||_{L^\infty} $ } &
{ $\| e^{0} \|_{L^2(L^2)} $ } &
{ $\| e^{1} \|_{L^2(L^2)} $ } &
{ $||| E |||_{L^2} $ } \\
\cmidrule(lr{.5em}){1-1}
\cmidrule(lr{.5em}){2-2}
\cmidrule(lr{.5em}){3-3}
\cmidrule(lr{.5em}){4-4}
\cmidrule(lr{.5em}){5-5}
\cmidrule(lr{.5em}){6-6}
\cmidrule(lr{.5em}){7-7}
\cmidrule(lr{.5em}){8-8}
$\tau_0/2^0$ & $h_0/2^0$ & 2.834e-02 & 2.862e-01 & 6.122e-01 & 2.099e-02 & 2.234e-01 & 4.808e-01 \\
$\tau_0/2^1$ & $h_0/2^1$ & 1.383e-03 & 1.755e-02 & 5.343e-02 & 9.773e-04 & 1.186e-02 & 3.989e-02 \\
$\tau_0/2^2$ & $h_0/2^2$ & 9.261e-05 & 1.075e-03 & 6.750e-03 & 6.064e-05 & 7.140e-04 & 4.835e-03 \\
$\tau_0/2^3$ & $h_0/2^3$ & 5.911e-06 & 6.690e-05 & 8.466e-04 & 3.812e-06 & 4.446e-05 & 6.005e-04 \\
$\tau_0/2^4$ & $h_0/2^4$ & 3.714e-07 & 4.186e-06 & 1.059e-04 & 2.387e-07 & 2.777e-06 & 7.495e-05 \\
$\tau_0/2^5$ & $h_0/2^5$ & 2.325e-08 & 2.616e-07 & 1.324e-05 & 1.492e-08 & 1.735e-07 & 9.364e-06 \\
\midrule 
\multicolumn{2}{c}{EOC} & 4.00 & 4.00 & 3.00 & 4.00 & 4.00 & 3.00\\
\bottomrule 
\end{tabular}
}
\end{table}

In our first we study the convergence behavior of the Galerkin--collocation approximation 
$U_{\tau,h}\in \left(X_\tau^3\big( V_h^{(3)}\big)\cap C^1\big(I;V_h^{(3)}\big)\right)^2$ 
of 
Problem \ref{Prob:DisLocalcGPC} for the prescribed solution 
\begin{equation}
\label{Eq:ExSol1}
u(\boldsymbol{x},t) := \sin(4 \pi t) \cdot \sin(2\pi x_1) \cdot \sin(2\pi x_2)
\end{equation}
of the wave problem \eqref{Eq:IBVP} on the space-time domain $\Omega \times I = (0,1)^2 
\times (0,1)$. For the piecewise polynomial order in space and time of the finite element 
approach the choice $k=3$ and $r=3$ is thus made; cf.\ \eqref{Eq:DefXk} and 
\eqref{Eq:DefVh}. Beyond the norms of $L^\infty(I;L^2(\Omega))$ and $L^2(I;L^2(\Omega))$ 
the convergence behavior is studied further with respect to the energy quantities
\begin{equation}
\label{Def:EN_1}
||| E |||_{L^\infty} = \max_{t \in \mathbb{I}}
( \| \nabla e^{0}(t) \|^2 + \| e^{1}(t) \|^2 )^{1/2} \quad \text{and} \quad 
||| E |||_{L^2} = \Big( \int_I
( \| \nabla e^{0} (t) \|^2 + \| e^{1}(t) \|^2)
\mathrm{d} t \Big)^{1/2}
\end{equation}
with $E(t) = U(t) - U_{\tau,h}(t)$. Throughout, the $L^\infty$-norms in time are computed 
on the discrete time grid
\[
\mathbb{I} = \{ t_n^j \::\: t_n^j = t_{n-1} + j \cdot k_n \cdot \tau_n,\,
k_n=0.001,\, j=0,\,\dots,999,\, n=1,\dots,N \} \cup \{ t_N \}\,.
\]
In the numerical experiments the domain $\Omega$ is decomposed into a sequence of 
successively refined meshes $\Omega_h^l$, with $l= 0,\ldots ,4$, of quadrilateral finite 
elements. On the coarsest level, we use a uniform decomposition of $\Omega$ 
into $4$ cells, corresponding to the mesh size $h_0=1/\sqrt{2}$, and of the time 
interval $I$ into $N=10$ subintervals which amounts to the time step size $\tau_0=0.1$. 
In the experiments the temporal and spatial mesh sizes are successively refined by a 
factor of two in each refinement step. 

In Table~\ref{Tab:1} we summarize the calculated results for this experiment. The 
experimental order of convergence (EOC) was calculated using the results from the two 
finest meshes. The 
numerical results of Table~\ref{Tab:1} nicely confirm our error estimates 
\eqref{Eq:ErrE_01} and 
\eqref{Eq:ErrE_02} by depicting the expected optimal fourth order rate of convergence in 
space and time. The third order convergence of the energy errors \eqref{Def:EN_1} is in 
agreement with the error estimates \eqref{Eq:ErrE_01a} and \eqref{Eq:ErrE_022}. Increasing 
the piecewise polynomial order in space to $r=4$ and thus considering an approximation 
$U_{\tau,h}\in \left(X_\tau^3 \big(V_h^{(4)}\big)\cap C^1\big(I;V_h^{(4)}\big)\right)^2$ 
in 
Problem \ref{Prob:DisLocalcGPC} leads a fourth order convergence behavior in time and 
space which is not shown here for the sake of limited space.

\subsection{Convergence test for
$\boldsymbol{\mathrm{cGP\text{-}C^1(4)}}$--$\boldsymbol{\mathrm{cG}(5)}$
and post-processing}

\begin{table}[htb!]
\centering
\caption{Error $E = \big(e^0, e^1\big) = U - U_{\tau,h}$ and error
$\widetilde{E} = \big(\tilde{e}^0, \tilde{e}^1\big) = U - \widetilde{U}_{\tau,h}$
of the post-processed solution $\widetilde{U}_{\tau,h}$ of Thm.~\ref{Thm:PP}, both with 
the
corresponding experimental orders of convergence (EOC),
for the solution $U = (u, \dt u)$ of \eqref{Eq:ExSol2} and the
Galerkin--collocation approximation $U_{\tau,h}\in \left(X_\tau^4
\big(V_h^{(5)}\big)\cap C^1\big(I;V_h^{(5)}\big)\right)^2$ of
Problem~\ref{Prob:DisLocalcGPC}.}
\label{Tab:2A}
{
\begin{tabular}{cccccccc}
\toprule
{$\tau$} & {$h$} &
{ $\| e^{0} \|_{L^\infty(L^2)} $ } &
{ $\| e^{1} \|_{L^\infty(L^2)} $ } &
{ $||| E |||_{L^\infty} $ } &
{ $\| e^{0} \|_{L^2(L^2)} $ } &
{ $\| e^{1} \|_{L^2(L^2)} $ } &
{ $||| E |||_{L^2} $ } \\
\cmidrule(lr{.5em}){1-1}
\cmidrule(lr{.5em}){2-2}
\cmidrule(lr{.5em}){3-3}
\cmidrule(lr{.5em}){4-4}
\cmidrule(lr{.5em}){5-5}
\cmidrule(lr{.5em}){6-6}
\cmidrule(lr{.5em}){7-7}
\cmidrule(lr{.5em}){8-8}
$\tau_0/2^0$ & $h_0$ & 8.457e-06 & 9.634e-05 & 9.637e-05 & 4.787e-06 & 5.392e-05 & 5.806e-05 \\
$\tau_0/2^1$ & $h_0$ & 2.497e-07 & 3.018e-06 & 3.022e-06 & 1.360e-07 & 1.654e-06 & 1.763e-06 \\
$\tau_0/2^2$ & $h_0$ & 7.608e-09 & 9.368e-08 & 9.372e-08 & 4.127e-09 & 5.141e-08 & 5.463e-08 \\
$\tau_0/2^3$ & $h_0$ & 2.353e-10 & 2.936e-09 & 2.936e-09 & 1.280e-10 & 1.604e-09 & 1.703e-09 \\
$\tau_0/2^4$ & $h_0$ & 7.323e-12 & 9.175e-11 & 9.175e-11 & 3.991e-12 & 5.012e-11 & 5.321e-11 \\
\midrule 
\multicolumn{2}{c}{EOC} & 5.01 & 5.00 & 5.00 & 5.00 & 5.00 & 5.00\\
\midrule\midrule
{$\tau$} & {$h$} &
{ $\| \tilde{e}^{0} \|_{L^\infty(L^2)} $ } &
{ $\| \tilde{e}^{1} \|_{L^\infty(L^2)} $ } &
{ $||| \widetilde{E} |||_{L^\infty} $ } &
{ $\| \tilde{e}^{0} \|_{L^2(L^2)} $ } &
{ $\| \tilde{e}^{1} \|_{L^2(L^2)} $ } &
{ $||| \widetilde{E} |||_{L^2} $ } \\
\cmidrule(lr{.5em}){1-1}
\cmidrule(lr{.5em}){2-2}
\cmidrule(lr{.5em}){3-3}
\cmidrule(lr{.5em}){4-4}
\cmidrule(lr{.5em}){5-5}
\cmidrule(lr{.5em}){6-6}
\cmidrule(lr{.5em}){7-7}
\cmidrule(lr{.5em}){8-8}
$\tau_0/2^0$ & $h_0$ & 2.906e-06 & 1.711e-05 & 1.791e-05 & 1.936e-06 & 1.519e-05 & 
1.764e-05 \\
$\tau_0/2^1$ & $h_0$ & 4.717e-08 & 2.802e-07 & 2.841e-07 & 3.150e-08 & 2.418e-07 & 
2.824e-07 \\
$\tau_0/2^2$ & $h_0$ & 7.513e-10 & 4.507e-09 & 4.537e-09 & 4.972e-10 & 3.797e-09 & 
4.440e-09 \\
$\tau_0/2^3$ & $h_0$ & 1.180e-11 & 7.085e-11 & 7.133e-11 & 7.788e-12 & 5.940e-11 & 
6.949e-11 \\
$\tau_0/2^4$ & $h_0$ & 1.851e-13 & 1.113e-12 & 1.120e-12 & 1.216e-13 & 9.282e-13 & 
1.086e-12 \\
\midrule 
\multicolumn{2}{c}{EOC} & 6.00 & 6.00 & 6.00 & 6.00 & 6.00 & 6.00\\
\bottomrule 
\end{tabular}
}
\end{table}


In the second numerical experiment we study the Galerkin--collocation scheme of 
Problem~\ref{Prob:DisLocalcGPC} for $k=4$ to obtain a fully discrete solution
$U_{\tau,h}\in \left(X_\tau^4\big( V_h^{(5)}\big)\cap C^1\big(I;V_h^{(5)}\big)\right)^2$.
In addition, we will apply the post-processing considered in Sect.~\ref{Sec:PP} and obtain a
solution $\widetilde{U}_{\tau,h}$ belonging to 
$\left(X_\tau^5\big( V_h^{(5)}\big)\cap C^2\big(I;V_h^{(5)}\big)\right)^2$.
The numerical study is done for the prescribed solution 
\begin{equation}
\label{Eq:ExSol2}
u(\boldsymbol{x},t) := \sin(4 \pi t) \, x_1  (1-x_1)  x_2  (1-x_2)
\end{equation}
of the wave problem \eqref{Eq:IBVP} on the space-time domain $\Omega \times I = (0,1)^2 
\times (0,1)$. For the piecewise polynomial order in space and time of the finite 
element approach the choice $k=4$ and $r=5$ is thus made. Since this work focusses
on the temporal discretization, the polynomial degree $r$ in space is chosen such that
the spatial approximation becomes exact. Hence, the convergence behavior in time
can be illustrated on a fixed spatial grid that consists of $4\times 4$ congruent
squares leading to $h_0=0.25\sqrt{2}$. The largest time length is $\tau_0=0.1$.

In Table~\ref{Tab:2A} we summarize the calculated results for this experiment. The
experimental order of convergence (EOC) was determined from the results on the two
finest meshes. The numerical results of Table~\ref{Tab:2A} nicely confirm the fifth
order rate of convergence of the $\mathrm{cGP\text{-}C^{1}(4)}$ time discretization. The 
application of the post-processing presented in Theorem~\ref{Thm:PP} increased all 
convergence rates from $5$ to $6$. This order can at most be expected for a polynomial 
approximation in time with piecewise polynomials of fifth order. By means of 
Theorem~\ref{Thm:PP}, Table~\ref{Tab:2A} thus underlines the optimal order 
approximation properties of the $\mathrm{cGP\text{-}C^{2}(5)}$ member of the
family of Galerkin--collocation schemes of Problem~\ref{Prob:PP}.

If the ${\mathrm{cGP\text{-}C^2(5)}}$--${\mathrm{cG}(5)}$ method of  
Problem~\ref{Prob:PP} is directly applied for the computation, instead of using the 
post-processing of Thm.~\ref{Thm:PP}, then exactly the same errors as shown in 
Table~\ref{Tab:2A} for $\widetilde{E}$ are obtained. However, using the post-processing 
has certain computational advantages. Since the
cGP-C$^1$($k$) approach leads to system matrices of simpler block structure
compared to the cGP-C$^2$($k$+1) method, the construction of efficient
preconditioners simplifies; cf.\ \cite{AB19} for details.  

\section{Summary}

In this work we presented a family of space-time finite element methods for wave 
problems. The schemes combine the concepts of collocation methods and Galerkin 
approximation. Continuously differentiable in time fully discrete solutions were 
obtained. An optimal order error analysis was provided for this class of methods. By an 
direct extension of the construction principle a further class of schemes with twice 
continuously differentiable in time discrete solutions was presented. A
theorem 
regarding the connection of the two classes of schemes to each other by means of a 
post-processing was given. The proven error estimates and the expected convergence rates 
for the second class of schemes were illustrated by numerical experiments. The 
construction of the methods can be transferred to further classes of non-stationary 
partial differential equations. The approach offers large potential for the approximation 
of multi-physics problems in that the coefficient functions of the subproblems are given by 
the solutions and their time derivatives of coupled further subproblems. In addition, the 
presented post-processing can nicely be exploited for a posteriori error control and 
adaptive refinement of the temporal mesh.

\end{document}